\documentclass[12pt,a4paper]{article}

\usepackage[dvips]{graphicx}
\usepackage{latexsym}
\usepackage{amssymb} 
\usepackage{amsmath}
\usepackage{amscd}
\newtheorem{Th}{Theorem}[section]
\newtheorem{Co}[Th]{Corollary}

\newtheorem{Lem}[Th]{Lemma}

\newtheorem{Rem}[Th]{Remark}
\newtheorem{Def}{ Definition}[section]
\newtheorem{Pro}[Th]{Proposition}
\newcommand{\demo}{\par\noindent{\it Proof. \/}\ }
\newcommand{\enD}{\hfill $\Box$ \vspace{3truemm}\par}
\newcommand{\bx}{\mbox{\boldmath $x$}}
\newcommand{\bX}{\mbox{\boldmath $X$}}
\newcommand{\be}{\mbox{\boldmath $e$}}

\newcommand{\bv}{\mbox{\boldmath $v$}}
\newcommand{\sbv}{\mbox{\scriptsize \boldmath$v$}}
\newcommand{\sbxi}{\mbox{\scriptsize \boldmath$\xi$}}
\newcommand{\by}{\mbox{\boldmath $y$}}
\newcommand{\bo}{\mbox{\boldmath $0$}}

\newcommand{\bw}{\mbox{\boldmath $w$}}
\newcommand{\bN}{\mbox{\boldmath $N$}}
\newcommand{\bxi}{\mbox{\boldmath $\xi$}}

\newcommand{\bgamma}{\mbox{\boldmath $\gamma$}}
\newcommand{\bsigma}{\mbox{\boldmath $\sigma$}}
\newcommand{\bn}{\mbox{\boldmath $n$}}

\newcommand{\R}{{\mathbb R}}
\newcommand{\lon}{\longrightarrow}

\begin{document}
\title{\bf Total lightcone curvatures of spacelike submanifolds in Lorentz-Minkowski space}
\author{ Shyuichi IZUMIYA \thanks{Work partially  supported  by  } 
}

\date{\today}
\maketitle

\begin{abstract}
We introduce the totally absolute lightcone curvature for 
a spacelike submanifold with general codimension
and investigate global properties of this curvature.
One of the consequences is that the Chern-Lashof type inequality holds.
Then the notion of lightlike tightness is naturally induced.
\end{abstract}
\renewcommand{\thefootnote}{\fnsymbol{footnote}}
\footnote[0]{2010 Mathematics Subject classification:53C40, 53C42, 53C80 }
\footnote[0]{Key Words and Phrases: lightcone Gauss map, normalized lightcone
Killing-Lipschitz curvature, Chern-Lashof type
 Theorem, lightlike tight spacelike immersion}
\section{Introduction}
\par
In this paper we consider global properties of spacelike submanifolds in
Lorentz-Minkowski space.
The study of the extrinsic differential geometry of submanifolds in Lorentz-Minkowski space
is of interest in the special relativity theory.
Moreover, it is a natural generalization of the extrinsic geometry of submanifolds in Euclidean space.
In \cite{IzCar07} the case of codimension two spacelike submanifolds has been considered.
The normalized lightcone Gauss map was introduced which plays the similar role to
the Gauss map of a hypersurface in the Euclidean space.
For example, the Gauss-Bonnet type theorem holds for the corresponding Gauss-Kronecker curvature 
(cf., \cite[Theorem 6.5]{IzCar07}).
Moreover, we recently discovered a new geometry on the hyperbolic space which is different from
the Gauss-Bolyai-Lobachevskii geometry (i.e., the hyperbolic geometry) \cite{BIR10,BIR11,Izu,IzCar06}.
We call this new geometry the {\it horospherical geometry}.
The horospherical Gauss map  (or, the hyperbolic Gauss map) is one of the key notions in
the horospherical geometry. We also showed that the Gauss-Bonnet type theorem holds
for the horospherical Gauss-Kronecker curvature\cite{IzCar06}.
The notion of normalized lightcone Gauss maps unifies both the notion of Gauss maps in 
the Euclidean space
and the notion of horospherical Gauss maps in the hyperbolic space.
\par
In this paper we generalize the normalized lightcone Gauss map and the corresponding curvatures for general 
spacelike submanifolds in 
Lorentz-Minkowski space.
If we try to develop this theory as a direct analogy to the Euclidean case, there exist several problems.
The main problem is that the fiber of the unit normal bundle of a spacelike submanifold is a union of
the pseudo-spheres which is not only non-compact but also non-connected.
So, we can not integrate the curvatures along the fiber at each point.
Therefore, we cannot define the Lipschitz-Killing curvature analogous to
the Euclidean case directly \cite{CL57}.
In order to avoid this problem, we arbitrary choose a future directed unit normal vector field
along the submanifold
and consider the pseudo-orthonormal space of this timelike vector on each fiber of the 
normal bundle. Then we obtain a spacelike codimension two unit normal sphere bundle in the normal bundle over 
the submanifold whose fiber is the Euclidean sphere.
As a consequence, we define the normalized lightcone Lipschitz-Killing curvature
and the total absolute lightcone curvature at each point.
We remark that the values of these curvatures are not invariant under the Lorentzian motions.
However, the flatness with respect to the curvature is an invariant property.
We can show that the total absolute lightcone curvature is independent of the choice of
the unit future directed timelike normal vector field (cf., Lemma 6.2).
Although these curvatures are not Lorentzian invariant, we show that the Chern-Lashof type inequality holds for this curvature (cf, \S 7).
In \S 8 we consider codimension two spacelike submanifolds.
In this case the situation is different from the higher codimensional case.
We have two different normalized lightcone Lipschitz-Killing (i.e., Gauss-Kronecker) curvatures at each point.
The corresponding total absolute normalized lightcone Lipschitz-Killing (i.e.,Gauss-Kronecker) curvatures are also different (cf., the remark 
after Theorem 8.3).
However, we also have the Chern-Lashof type inequality for each total absolute Lipschitz-Killing (i.e., Gauss-Kronecker) curvature.
Moreover, we consider the Willmore type integral (cf., \cite[Theorem 7.2.2]{Will}) of the lightcone mean curvature for spacelike surface in Lorentz-Minkowski $4$-space. 
Finally, we introduce the notion of the lightlike tightness which characterize
the minimal value of the total absolute lightcone curvature.
As a special case, we have the horo-spherical Chern-Lashof type inequality and horo-tight immersions
in the hyperbolic space
 \cite{BIR10,BIR11,SolTeu}.
Motivated by those arguments, we can introduce the notion of several kinds of tightness and tautness
depending on the causal characters which will be one of the subjects of a future program of the research.

\section{Basic concepts in Lorentz-Minkowski space}
\par
We introduce in this section some basic notions on Lorentz-Minkowski
$n+1$-space. For basic
concepts and properties, see \cite{Oneil}.
\par
Let $\R^{n+1}=\{(x_0,x_1,\dots  ,x_n)\ |\ x_i \in \R\ (i=0,1,\dots , n)\ \}$
be an $n+1$-dimensional cartesian space. For any
$\bx =(x_0,x_1,\dots ,x_n),\ \by =(y_0,y_1,\dots ,y_n)\in \R^{n+1},$
the {\it pseudo scalar product } of $\bx $ and $\by $ is defined by
$$
\langle \bx ,\by\rangle =-x_0y_0+\sum _{i=1}^n x_iy_i.
$$
We call $(\R^{n+1} ,\langle ,\rangle )$  {\it Lorentz-Minkowski
$n+1$-space} (or, simply {\it Minkowski $n+1$-space\/}. We write $\R^{n+1}_1$ instead of $(\R^{n+1} ,\langle
,\rangle )$. We say that a non-zero vector $\bx \in \R^{n+1}_1$ is
{\it spacelike, lightlike or timelike} if $\langle \bx,\bx \rangle
>0,$ $\langle \bx,\bx \rangle =0$ or  $\langle \bx,\bx \rangle <0$
respectively. The norm of the vector $\bx \in \R^{n+1}_1$ is defined
to be $\|\bx\|=\sqrt{|\langle \bx,\bx\rangle |}.$ We have the canonical
projection $\pi :\R^{n+1}_1\lon \R^n$ defined by $\pi (x_0,x_1,\dots
,x_n)=(x_1,\dots ,x_n).$ Here we identify $\{\bo\}\times \R^n$ with
$\R^n$ and it is considered as Euclidean $n$-space whose scalar
product is induced from the pseudo scalar product $\langle ,\rangle
.$ For a vector $\bv \in \R^{n+1}_1$ and a real number $c,$ we
define a {\it hyperplane with pseudo normal\/} $\bv$ by
$$
HP(\bv ,c)=\{\bx\in \R^{n+1}_1\ |\ \langle \bx ,\bv \rangle = c\ \}.
$$
We call $HP(\bv ,c)$ a {\it spacelike hyperplane\/}, a {\it timelike hyperplane\/}
or a {\it lightlike hyperplane\/}  if $\bv $ is timelike, spacelike or lightlike respectively.
\par
We now define {\it Hyperbolic $n$-space} by
$$
\mathbb{H}^n(-1)=\{\bx\in \R^{n+1}_1 | \langle \bx ,\bx\rangle =-1
\}
$$
and {\it de Sitter $n$-space} by
$$
\mathbb{S}^n_1=\{\bx\in \R^{n+1}_1 | \langle \bx ,\bx\rangle =1\ \}.
$$

We define
$$
\mathbb{LC}^*=\{\bx=(x_0,x_1,\dots ,x_n)\in \R^{n+1}_1  \ | x_0\not= 0,\
\langle\bx ,\bx \rangle =0\}
$$
and we call it {\it the {\rm (}open{\rm )} lightcone} at the origin.

If $\bx =(x_0,x_1,\dots ,x_2)$ is a non-zero lightlike vector, then
$x_0\not= 0.$ Therefore we have
\[
\widetilde{\bx}=\left(1,\frac{x_1}{x_0},\dots
,\frac{x_n}{x_0}\right)\in \mathbb{S}^{n-1}_+=\{\bx =(x_0,x_1,\dots ,x_n)\ |\
\langle \bx ,\bx \rangle =0,\ x_0=1\}.
\]
We call $\mathbb{S}^{n-1}_+$ the {\it lightcone unit $n-1$-sphere.}
\par
\par
For any $\bx_1,\bx_2,\dots ,\bx_n \in \R^{n+1}_1,$
we define a vector $\bx_1\wedge\bx_2\wedge\dots \wedge\bx_n$
by
\[
\bx_1\wedge\bx_2\wedge\dots \wedge\bx_n=
\left|
\begin{array}{cccc}
-\be_{0}&\be_{1}&\cdots &\be_{n}\\
x^1_{0}&x^1_{1}&\cdots &x^1_{n}\\
x^2_{0}&x^2_{1}&\cdots &x^2_{n}\\
\vdots &\vdots &\cdots &\vdots \\
x^n_{0}&x^n_{1}&\cdots &x^n_{n}
\end{array}
\right| ,
\]
where $\be_{0},\be_{1},\dots ,\be_{n}$ is the canonical basis of $\R^{n+1}_1$
and $\bx_i=(x_0^i,x_1^i,\dots ,x_n^i).$
We can easily check that
$$
\langle \bx,\bx_1\wedge\bx_2\wedge\dots \wedge\bx_n\rangle ={\rm det}(\bx,\bx_1,\dots ,\bx_n),
$$
so that
$\bx_1\wedge\bx_2\wedge\dots \wedge\bx_n$ is pseudo orthogonal to  any $\bx _i$ $(i=1,\dots ,n).$

\section{Differential geometry on spacelike
submanifolds}
\par
In this section we introduce the basic geometrical framework for the 
study of spacelike submanifolds in Minkowski $n+1$-space analogous to
the case of codimension two in \cite{IzCar07}.
Let $\R^{n+1}_1$ be an oriented and time-oriented
space. We choose $\be _0=(1,0,\dots ,0)$ as the future timelike
vector field. 
Let $\bX :U\lon \R^{n+1}_1$ be a spacelike embedding of codimension $k,$
where $U\subset \R^s$ ($s+k=n+1$) is an open subset. We also write $M=\bX(U)$ and identify $M$ and $U$ through the embedding $\bX.$
We say that $\bX$ is {\it spacelike} if the tangent space $T_p  M$
of $M$ at $p$ is a spacelike subspace (i.e., consists of spacelike
vectors) for any point $p\in M$.
For any $p=\bX(u)\in M\subset \R^{n+1}_1,$
we have
\[
T_pM=\langle \bX_{u_1}(u),\dots ,\bX_{u_s}(u)\rangle _\R.
\]
Let $N_p(M)$ be the pseudo-normal space of $M$ at $p$ in $\R^{n+1}_1.$
Since $T_pM$ is a spacelike subspace of $T_p\R^{n+1}_1,$
$N_p(M)$ is a $k$-dimensional Lorentzian subspace of $T_p\R^{n+1}_1$
(cf.,\cite{Oneil}).
On the pseudo-normal space $N_p(M),$ we have two kinds of pseudo spheres:
\begin{eqnarray*}
N_p(M;-1)& = & \{\bv\in N_p(M)\ |\ \langle \bv,\bv\rangle =-1\ \} \\
N_p(M;1)&= & \{\bv\in N_p(M)\ |\ \langle \bv,\bv\rangle =1\ \},
\end{eqnarray*}
so that we have two unit spherical normal bundles over $M$:
\[
N(M;-1)=\bigcup _{p\in M} N_p(M;-1)\ \mbox{and}\  N(M;1)=\bigcup _{p\in M} N_p(M;1).
\]
Then we have the Whitney sum decomposition
\[
T\R^{n+1}_1|_{M}=TM\oplus N(M).
\]
Since $M=\bX(U)$ is spacelike, $\be _0$ is a transversal future directed timelike vector field along $M$.
For any $\bv\in T_p\R^{n+1}_1|_{M},$ we have $\bv=\bv_1+\bv_2,$ where $\bv_1\in T_pM$ and $\bv_2\in N_p(M).$
If $\bv$ is timelike, then $\bv_2$ is timelike. Let $\pi _{N(M)}:T\R^{n+1}_1|_{M}\lon N(M)$ be the canonical projection.
Then $\pi_{N(M)}(\be_0)$ is a future directed timelike normal vector field along $M.$
So we always have a future directed unit timelike normal vector field along $M$ (even globally).
We now arbitrarily choose a future directed unit timelike normal vector field  $\bn ^T(u)\in N_p(M;-1),$ where $p=\bX(u).$
Therefore we have the pseudo-orthonormal complement
$(\langle \bn ^T(u)\rangle _\R)^\perp$ in $N_p(M)$
which is a $k-1$-dimensional subspace of $N_p(M).$
We can also choose a pseudo-normal section 
$\bn ^S(u)\in (\langle \bn ^T(u)\rangle _\R)^\perp\cap N(M;1)$ 
at least locally, then we have
$\langle \bn ^S,\bn ^S\rangle =1$ and $\langle \bn ^S,\bn ^T\rangle =0.$
We define a $k-1$-dimensional spacelike unit sphere in $N_p(M)$ by
\[
N_1(M)_p[\bn ^T]=\{\bxi \in N_p(M;1)\ |\ \langle \bxi, \bn ^T (p)\rangle =0\ \}.
\]
Then we have a {\it spacelike unit $k-1$-spherical bundle over $M$ with respect to $\bn ^T$} defined by
\[
N_1(M)[\bn ^T]=\bigcup _{p\in M} N_1(M)_p[\bn ^T].
\]
Since we have
$T_{(p,\xi)}N_1(M)[\bn^T]=T_pM\times T_\xi N_1(M)_p[\bn ^T],$
we have the canonical Riemannian metric on $N_1(M)[\bn^T].$
We denote the Riemannian metric on $N_1(M)[\bn^T]$ by $(G_{ij}(p,\bxi))_{1\leqslant i,j\leqslant n-1}.$
\par
For any future directed unit normal $\bn^T$ along $M,$ we arbitrary choose the unit
spacelike normal vector field $\bn^S$ with $\bn^S(u)\in N_1(M)_p[\bn^T]$, where $p=\bX(u).$
We call $(\bn^T,\bn^S)$ a {\it future directed pair\/} along $M.$
Clearly, the vectors
$\bn^T (u)\pm \bn^S(u)$ are lightlike. Here we choose $\bn ^T+\bn^S$
as a lightlike normal vector field along $M.$ 
\begin{Def}{\rm
We define a mapping
\[
\mathbb{LG}(\bn^T,\bn^S):U\lon \mathbb{LC}^*
\]
by $\mathbb{LG}(\bn^T,\bn^S)(u)=\bn^T(u)+\bn^S(u).$
We call it the {\it lightcone Gauss image} of $M=\bX(U)$ with respect to
$(\bn^T,\bn^S).$
We also define a mapping
\[
\widetilde{\mathbb{LG}}(\bn^T,\bn^S): U\lon \mathbb{S}^{n-1}_+
\]
by $\widetilde{\mathbb{LG}}(\bn^T,\bn^S)(u)=\widetilde{\bn^T(u)+\bn^S(u)}$
which is called the {\it lightcone Gauss map} of $M=\bX(U)$ with respect to
$(\bn^T,\bn^S).$ 
}
\end{Def}
Under the identification of $M$ and $U$ through $\bX,$ we have the
linear mapping provided by the derivative of the lightcone Gauss image $\mathbb{LG}(\bn^T,\bn^S)$ at each point $p\in M$,
\[
d_p\mathbb{LG}(\bn^T,\bn^S):T_pM\lon T_p\R^{n+1}_1= T_pM\oplus N_p(M).
\]
Consider the orthogonal projections $\pi ^{\tau}:T_pM\oplus
N_p(M)\rightarrow T_p(M)$ and $\pi ^{\nu}:T_p(M)\oplus N_p(M)\rightarrow
N_p(M).$ We define
\[
d_p\mathbb{LG}(\bn^T,\bn^S)^{\tau}=\pi ^{\tau}\circ d_p(\bn^T+\bn^S)
\]
 and
\[
d_p\mathbb{LG}(\bn^T,\bn^S)^{\nu}=\pi ^{\nu}\circ d_p(\bn^T+\bn^S).
\]
\begin{Def}{\rm  We respectively call the
linear transformations $S_{p}(\bn^T,\bn^S)=- d_p\mathbb{LG}(\bn^T,\bn^S)^{\tau}$ and
$d_p\mathbb{LG}(\bn^T,\bn^S)^{\nu}$ of $T_{p}M$, the {\it $(\bn^T,\bn^S)$-shape
operator} of $M=\bX (U)$ at $p=\bX (u)$ and the {\it normal
connection with respect to $(\bn^T,\bn^S)$} of $M=\bX (U)$ at $p=\bX
(u).$
The eigenvalues of $S_{p}(\bn^T,\bn^S)$, denoted by $\{\kappa
_{i}(\bn^T,\bn^S)(p)\}_ {i=1}^s,$ are called the {\it lightcone
principal curvatures  with respect to $(\bn^T,\bn^S) $\/} at $p=\bX(u)$.
Then the {\it lightcone Gauss-Kronecker curvature with respect to
$(\bn^T,\bn^S)$\/} at $p=\bX (u)$ is defined by
\[
K_\ell(\bn^T,\bn^S)(p)={\rm det} S_{p}(\bn^T,\bn^S).
\]
We say that a point $p=\bX (u)$ is an {\it $(\bn^T,\bn^S)$-umbilical
point} $S_{p}(\bn^T,\bn^S)=\kappa (\bn^T,\bn^S)(p) 1_{T_{p}M}$.  We say that $M=\bX (U)$ is {\it totally
$(\bn^T,\bn^S)$-umbilical} if all points on $M$ are
$(\bn^T,\bn^S)$-umbilical.
}
\end{Def}
\par
We deduce now the lightcone Weingarten formula. Since $\bX _{u_i}$
$(i=1,\dots s)$ are spacelike vectors, we have a Riemannian metric
(the {\it first fundamental form \/}) on $M=\bX (U)$
defined by $ds^2 =\sum _{i=1}^{s} g_{ij}du_idu_j$,  where
$g_{ij}(u) =\langle \bX _{u_i}(u ),\bX _{u_j}(u)\rangle$ for any
$u\in U.$ We also have a {\it lightcone second fundamental invariant
with respect to the normal vector field $(\bn^T,\bn ^S) $\/} defined
by $h _{ij}(\bn^T,\bn^S )(u)=\langle -(\bn^T +\bn^S)
_{u_i}(u),\bX_{u_j}(u)\rangle$ for any $u\in U.$
By the similar arguments to those in the proof of \cite[Proposition 3.2]{IzCar07}, we have 
the following proposition.
\begin{Pro}
We choose a pseudo-orthonormal frame $\{\bn^T,\bn^S_1,\dots ,\bn^S_{k-1}\}$ of $N(M)$ with $\bn^S_{k-1}=\bn^S.$ Then we have the following lightcone Weingarten formula with respect to $(\bn^T,\bn^S)${\rm :}
\vskip1.5pt
\par
{\rm (a)} $\mathbb{LG}(\bn^T,\bn^S)_{u_i}=\langle \bn ^S_{u_i},\bn ^T\rangle(\bn^T-\bn^S)+\sum _{\ell =1}^{k-2}\langle (\bn^T+\bn^S)_{u_i},\bn^S_\ell \rangle\bn^S_\ell -\sum_{j=1}^{s}
h_i^j(\bn^T,\bn^S )\bX _{u_j}$
\par
{\rm (b)} $
\pi ^{\tau}\circ \mathbb{LG}(\bn^T,\bn^S)_{u_i}=-\sum_{j=1}^{s}
h_i^j(\bn^T,\bn^S )\bX _{u_j}.
$
\smallskip
\par\noindent
Here $\displaystyle{\left(h_i^j(\bn^T,\bn^S )\right)=\left(h_{ik}(\bn^T,\bn^S)\right)\left(g^{kj}\right)}$
and $\displaystyle{\left( g^{kj}\right)=\left(g_{kj}\right)^{-1}}.$
\end{Pro}
\par
As a corollary of the above proposition, we have an explicit
expression of the lightcone curvature in terms of the Riemannian
metric and the lightcone second fundamental invariant.
\begin{Co}
Under the same notations as in the above proposition, the lightcone
Gauss-Kronecker curvature relative to $(\bn^T,\bn^S) $ is given
by
$$
K_\ell (\bn^T,\bn^S )=\frac{\displaystyle{{\rm det}\left(h_{ij}(\bn^T,\bn^S )\right)}}
{\displaystyle{{\rm det}\left(g_{ij}\right)}}.
$$
\end{Co}
\par
Since $\langle -(\bn^T +\bn^S )(u),\bX _{u_j}(u)\rangle =0,$ we have
$h_{ij}(\bn ^T,\bn^S)(u)=\langle \bn^T (u)+\bn^S (u),\bX
_{u_iu_j}(u)\rangle.$ Therefore the lightcone second fundamental
invariant at a point $p_0=\bX (u_0)$ depends only on the values 
$\bn^T (u_0)+\bn^S (u_0)$ and $\bX _{u_iu_j}(u_0)$, respectively.
Thus, the lightcone curvatures also depend only on
$\bn^T (u_0)+\bn^S (u_0)$, $\bX_{u_i}(u_0)$  and $\bX
_{u_iu_j}(u_0)$, independent of the derivation of the vector fields 
$\bn^T$ and $\bn^S .$ We write $\kappa _i(\bn^T_0,\bn^S_0)(p_0)$ $(i=1,\dots ,s)$
and $K_\ell (\bn
^T_0,\bn^S_0)(u_0)$ as the lightcone curvatures at $p_0=\bX (u_0)$
with respect to $(\bn ^T_0,\bn^S_0)=(\bn^T (u_0),\bn^S(u_0)).$ We
might also say that a point $p_0=\bX (u_0)$ is 
{\it $(\bn^T_0,\bn^S_0)$-umbilical\/} because the lightcone $(\bn^T,\bn^S) $-shape
operator at $p_0$ depends only on the normal vectors $(\bn
^T_0,\bn^S_0).$
So we denote that $h_{ij}(\bn^T,\bxi)(u_0)=h_{ij}(\bn^T,\bn^S)(u_0)$ and 
 $K_\ell(\bn^T,\bxi)(p_0)=K_\ell(\bn^T_0,\bn^S_0)(p_0)$,
where $\bxi =\bn^S(u_0)$ for some local extension $\bn^T(u)$ of $\bxi.$
Analogously, we say that a point $p_0=\bX (u_0)$ is an {\it $(\bn
^T_0,\bn^S_0)$-parabolic point \/} of $\bX :U\lon \R^{n+1}_1$ if
$K_\ell (\bn ^T_0,\bn^S_0)(u_0)=0.$ And we say that a point $p_0=\bX
(u_0)$ is a {\it $(\bn ^T_0,\bn^S_0)$-flat point \/} if it is an
$(\bn ^T_0,\bn^S_0)$-umbilical point and $K_\ell(\bn^T
_0,\bn^S_0)(u_0)=0.$
\par
On the other hand, the lightcone Gauss map $\widetilde{\mathbb{LG}}(\bn^T,\bn^S)$ with respect to
$(\bn^T,\bn^S)$ also induces a linear mapping $d_p\widetilde{\mathbb{LG}}(\bn^T,\bn^S):T_pM\lon T_p\R^{n+1}_1$ under
the identification of $U$ and $M,$ where $p=\bX(u).$
We have the following proposition.
\begin{Pro}
Under the above notations, we have the following normalized lightcone Weingarten formula with respect to $(\bn^T,\bn^S)$:
\[
\pi^{\tau}\circ \widetilde{\mathbb{LG}}(\bn^T,\bn^S)_{u_i}=-\sum_{j=1}^s \frac{1}{\ell _0(u)}h^j_i(\bn^T,\bn^S)\bX_{u_j},
\]
where $\mathbb{LG}(\bn^T,\bn^S)(u)=(\ell _0(u),\ell _1(u),\dots ,\ell _n(u)).$
\end{Pro}
\demo
By definition, we have 
$\ell_0\widetilde{\mathbb{LG}}(\bn^T,\bn^S)=\mathbb{LG}(\bn^T,\bn^S).$ It follows that
$\ell _0\widetilde{\mathbb{LG}}(\bn^T,\bn^S)_{u_i}=\mathbb{LG}(\bn^T,\bn^S)_{u_i}-\ell _{0u_i}\widetilde{\mathbb{LG}}(\bn^T,\bn^S).$
Since $\widetilde{\mathbb{LG}}(\bn^T,\bn^S)(u)\in N_p(M),$ we have
\[
\pi^{\tau}\circ\widetilde{\mathbb{LG}}(\bn^T,\bn^S)_{u_i}=\frac{1}{\ell _0}\pi^{\tau}\circ\mathbb{LG}(\bn^T,\bn^S)_{u_i}.
\]
By the lightcone Weingarten formula with respect to $(\bn^T,\bn^S)$ (Proposition 3.1), we have the desired formula.
\enD
\begin{Def}{\rm 
We call the linear transformation $\widetilde{S}(\bn^T,\bn^S)_p=-\pi^{\tau}\circ d_p\widetilde{\mathbb{LG}}(\bn^T,\bn^S)$ the {\it normalized
lightcone shape operator} of $M$ at $p$ with respect to $(\bn^T,\bn^S).$
The eigenvalues $\{\widetilde{\kappa} _i(\bn^T,\bn^S)(p)\}_{i=1}^s$ of 
$\widetilde{S}(\bn^T,\bn^S)_p$ are 
called the {\it normalized
lightcone principal curvatures}.
By the above proposition, we have
$\widetilde{\kappa} _i(\bn^T,\bn^S)(p)=(1/\ell _0(u))\kappa _i(\bn^T,\bn^S)(p).$
The {\it normalized Gauss-Kronecker curvature} of $M$ with respect to $(\bn^T,\bn^S)$ is defined to be
$$
\widetilde{K}_\ell (\bn^T,\bn^S)(u)={\rm det}\, \widetilde{S}(\bn^T,\bn^S)_p.
$$
}
\end{Def}
Then we have the following relation between the normalized lightcone Gauss-Kronecker curvature and the lightcone
Gauss-Kronecker curvature:
\[
\widetilde{K}_\ell(\bn^T,\bn^S)(u)=\left(\frac{1}{\ell_0(u)}\right)^sK_\ell (\bn^T,\bn^S)(u).
\]
\par
On the other hand, we consider a submanifold $\Delta =\{(\bv,\bw)\ |\ \langle \bv,\bw\rangle =0\ \}\subset H^n_+(-1)\times \mathbb{S}^n_1$ and the canonical projection $\bar{\pi}:\Delta \lon H^n_+(-1).$
It is well known that $\Delta$ can be identified with the unit tangent bundle $S(TH^n_+(-1))$ over
$H^n_+(-1).$
We define a function $\mathcal{N}_h:\Delta \lon \R$ by
$\mathcal{N}_h(\bv,\bw)=1/(v_0+w_0),$
where $\bv=(v_0,v_1,\dots,v_n),\bw=(w_0,w_1,\dots ,w_n).$
Then we have
\[
\mathcal{N}_h(\bn^T(u),\bn^S(u))=\frac{1}{\ell _0(u)}.
\]
Therefore we can rewrite the above formula as follows:
\[
\widetilde{K}_\ell(\bn^T,\bn^S)(u)=\mathcal{N}_h(\bn^T(u),\bn^S(u))^sK_\ell (\bn^T,\bn^S)(u).
\]

By definition, $p_0=\bX(u_0)$ is the $(\bn^T_0,\bn^S_0)$-umbilical point if and only if
$\widetilde{S}(\bn^T,\bn^S)_{p_0}=\widetilde{\kappa} _i(\bn^T,\bn^S)(p)1_{T_{p_0}M}.$
We have the following proposition.
\begin{Pro}
For a future directed unit normal vector field $\bn^T$ along $M=\bX(U),$
the following two conditions are equivalent{\rm :}
\par\noindent
{\rm (1)} There exists a spacelike unit normal vector field $\bn^S$ along $M=\bX(U)$ such that the normalized lightcone Gauss map $\widetilde{\mathbb{LG}}(\bn^T,\bn^S)$ of $M=\bX(U)$ is constant
\par\noindent
{\rm (2)} There exists $\bv\in \mathbb{S}^{n-1}_+$ and a real number $c$ such that $M\subset HP(\bv,c).$
\par
Suppose that the above condition holds. Then 
\par\noindent
{\rm (3)} $M=\bX(U)$ is totally $(\bn^T,\bn^S)$-flat.
\end{Pro}
\demo
Suppose that the normalized lightcone Gauss Map $\widetilde{\mathbb{LG}}(\bn^T,\bn^S)(u)=\bv$ is constant.
We consider a function $F:U\lon \R$ defined by $F(u)=\langle \bX(u),\bv\rangle.$
By definition, we have
\[
\frac{\partial F}{\partial u_i}(u)=\langle \bX_{u_i}(u),\bv\rangle
=\langle \bX_{u_i}(u),\widetilde{\mathbb{LG}}(\bn^T,\bn^S)(u)\rangle ,
\]
for any $i=1,\dots ,s.$ Therefore, $F(u)=\langle \bX(u),\bv\rangle =c$ is constant.
It follows that $M\subset HP(\bv, c)$ for $\bv\in \mathbb{S}^{n-1}_+.$
\par
Suppose that $M$ is a subset of a lightlike hyperplane $H(\bv,c)$ for $\bv\in \mathbb{S}^{n-1}_+.$
Since $M\subset HP(\bv ,c),$ we have $T_pM\subset H(\bv ,0)$.
If $\langle \bn^T(u),\bv\rangle=0,$ then $\bn^T(u)\in HP(\bv ,0).$
We remark that $HP(\bv ,0)$ does not contain timelike vectors.
This is a contradiction. So we have $\langle \bn^T(u),\bv\rangle \not=0.$
We now define a normal vector field along $M=\bX(U)$ by
\[
\bn^S(u)=\frac{-1}{\langle \bn^T(u),\bv\rangle}\bv-\bn^T(u).
\]
We can easily show that $\bn^S(u)\in  N_1(M)_{p}[\bn^T]$ for $p=\bX(u).$
Therefore $(\bn^T,\bn^S)$ is a future directed normal pair such that 
$\widetilde{\mathbb{LG}}(\bn^T,\bn^S)(u)=\bv.$
\par
On the other hand, by Proposition 3.3, if $\widetilde{\mathbb{LG}}(\bn^T,\bn^S)$ is constant, then
$(h^j_i(\bn^T,\bn^S)(u))=O$, so that $M=\bX(U)$ is lightcone $(\bn^T,\bn^S)$-flat.
\enD
\section{The normalized lightcone Lipschitz-Killing curvature}
In this section we define the lightcone Gauss map
of $N_1(M)[\bn^T]$ and investigate the geometric properties.
\begin{Def}
{\rm We define a map
\[
\widetilde{\mathbb{LG}}(\bn^T):N_1(M)[\bn^T]\lon \mathbb{S}^{n-1}_+
\]
by
$\widetilde{\mathbb{LG}}(\bn^T)(u,\bxi)=\widetilde{\bn^T(u)+\bxi},$
which we call the {\it lightcone Gauss map} of $N_1(M)[\bn^T].$
}
\end{Def}
The lightcone Gauss map leads us to a curvature similar to the codimension two case\cite{IzCar07}.
Let $T_{(p,\xi)}N_1(M)[\bn^T]$ be the tangent space of $N_1(M)[\bn^T]$ at $(p,\bxi).$
We have the canonical identification 
\[
T_{(p,\sbxi)}N_1(M)[\bn^T]=T_pM\oplus T_{\sbxi} \mathbb{S}^{k-2}\subset T_pM\oplus N_p(M)=T_p\R^{n+1}_1,
\]  
where $T_{\sbxi} \mathbb{S}^{k-2}\subset T_{\sbxi} N_p(M)\equiv N_p(M)$ and $p=\bX(u).$
Under this identification, we have
\[
T_p\R^{n+1}_1=T_pM\oplus N_p(M)=T_pM\oplus T_{\sbxi} \mathbb{S}^{k-2}\oplus \R^{k+1}=T_{(p,\sbxi)}N_1(M)[\bn^T]\oplus \R^{k+1}.
\]
Therefore, we can define the canonical projection 
\[
\Pi ^{\tau} :\widetilde{\mathbb{LG}}(\bn^T)^*T\R^{n+1}_1=TN_1(M)[\bn^T]\oplus \R^{k+1}
\lon TN_1(M)[\bn^T].
\]
 It follows that
we have a linear transformation
\[
\Pi^{\tau}_{\widetilde{\mathbb{LG}(n^T)(p,\xi)}}\circ d_{(p,\xi)}\widetilde{\mathbb{LG}}(\bn^T)
: T_{(p,\xi)}N_1(M)[\bn^T]\lon T_{(p,\xi)}N_1(M)[\bn^T].
\]
\begin{Def}{\rm 
The {\it normalized lightcone Lipschitz-Killing curvature} of $N_1(M)[\bn^T]$ at $(p,\bxi)$
is defined to be
\[
\widetilde{K}_\ell (\bn^T)(p,\bxi)=
\det \left(-\Pi^{\tau}_{\widetilde{\mathbb{LG}}(n^T) 
(p,\xi)}\circ d_{
(p,\xi)}\widetilde{\mathbb{LG}}(\bn^T)\right).
\]
\par
In order to investigate the lightcone Gauss map $\widetilde{\mathbb{LG}}(\bn^T)$
of $N_1(M)[\bn^T]$,
we define a map 
\[
\mathbb{LG}(\bn^T): N_1(M)[\bn^T]\lon \mathbb{LC}^*
\]
by
$\mathbb{LG}(\bn^T)(u,\bxi)=\bn^T(u)+\bxi,$
which  is called the {\it lightcone Gauss image} of $N_1(M)[\bn^T]$.
}
\end{Def}
We now write $\mathbb{LG}(\bn^T)(p,\bxi)=(\ell _0(p,\bxi),\ell _1(p,\bxi),\dots ,\ell _{n}(p,\bxi)).$ 
For any future directed timelike unit normal vector field $\bn^T$ along $M,$
there exists a pseudo-orthonormal frame $\{\bn^T,\bn^S_1,\dots ,\bn^S_{k-1}\}$
of $N(M)$ with $\bn^S_{k-1}(u_0)=\bxi$ and $p=\bX(u_0),$ so that we have a frame field
\[
\{\bX_{u_1},\dots, \bX_{u_s},\bn^T,\bn^S_1,\dots ,\bn^S_{k-1}\}
\]
of $\R^{n+1}_1$ along $M.$
We define an $\mathbb{S}^{k-2}$-family of spacelike unit normal vetor field 
\[
\bN(u,\mu)=\sum_{j=1}^{k-1}\mu _j\bn^S_j(u)
\in N(M;1)
\]
along $M$ for $\mu =(\mu_1,\dots ,\mu _{k-1})\in \mathbb{S}^{k-2}\subset \R^{k-1}.$ 
We also define a map
\[
\Psi : U\times \mathbb{S}^{k-2}\lon N_1(M)[\bn^T]
\]
by
$\Psi (u,\mu)=(\bX(u),\bN^S(u,\mu)),$ which gives a local parametrization
of $N_1(M)[\bn^T].$
Then we have $(p,\bxi)=(\bX(u_0),\bN^S(u_0,\mu_0)),$
where $\mu _0=(0,\dots ,0,1).$
It follows that 
$\mathbb{LG}(\bn^T)\circ \Psi (u,\mu)=\bn^T(u)+\bN^S(u,\mu).$
We now write that $\mathbb{LG}(\bn^T,\bN^S)(u,\mu)=\mathbb{LG}(\bn^T)\circ \Psi (u,\mu).$
We consider the local coordinate neighborhood of $\mathbb{S}^{k-2}$:
\[
U_{k-1}^+=\{(\mu_1,\dots ,\mu _{k-1})\in \mathbb{S}^{k-2}\ |\ \mu_{k-1}>0\ \}.
\]
Then we have $\mu_{k-1}=\sqrt{1-\sum_{j=1}^{k-2} \mu _j^2}.$
For $i=1,\dots, s, j=1,\dots k-2,$ we have the following calculation:
\begin{eqnarray*}
\frac{\partial \mathbb{LG}(\bn^T,\bN^S)}{\partial u_i}(u,\mu) &=& \bn^T_{u_i}(u)
+\sum_{\ell =1}^{k-1} \mu _{\ell}\bn ^S_{\ell ,u_i}(u), \\
\frac{\partial \mathbb{LG}(\bn^T,\bN^S)}{\partial \mu_j}(u,\mu) &=& \bn ^S_j(u)-
\frac{\mu _j}{\mu _{k-1}}\bn ^S_{k-1}(u).
\end{eqnarray*}
Therefore, we have 
\begin{eqnarray*}
\frac{\partial \mathbb{LG}(\bn^T,\bN^S)}{\partial u_i}(u_0,\mu _0) &=& \bn^T_{u_i}(u_0)
+\bn ^S_{k-1 ,u_i}(u_0)=(\bn^T+\bn^S_{k-1})_{u_i}(u_0), \\
\frac{\partial \mathbb{LG}(\bn^T,\bN^S)}{\partial \mu_j}(u_0,\mu _0) &=& \bn ^S_j(u_0).
\end{eqnarray*}
We now remark that $\{\bX_{u_1},\dots, \bX_{u_s},\bn^S_1,\dots ,\bn^S_{k-2}\}$ is 
a basis of $T_{(p,\xi)}N_1(M)[\bn^T]$ at $u=u_0.$
By Proposition 3.1, we have 
\begin{eqnarray*}
(\bn^T+\bn^S_{k-1})_{u_i}(u_0)&=& \langle \bn ^S_{u_i},\bn ^T\rangle(\bn^T-\bn^S_{k-1})(u_0) \\
&{}& +\sum _{\ell =1}^{k-2}\langle (\bn^T+\bn^S_{k-1})_{u_i},\bn^S_\ell \rangle\bn^S_\ell(u_0) 
-\sum_{j=1}^{s}
h_i^j(\bn^T,\bn^S )\bX _{u_j}(u_0).
\end{eqnarray*}
Since $\langle \bn^T-\bn^S_{k-1},\bX_{u_i}\rangle =
\langle \bn^T-\bn^S_{k-1},\bn^S_j\rangle =\langle \bn^S_\ell,\bX_{u_i}\rangle =0$
and $\langle \bn^S_j,\bn^S_{\ell}\rangle =\delta _{j\ell},$
we have
\begin{equation*}\label{1012061438}
 \begin{split}
& \det \left(-\Pi^t_{\mathbb{LG}(n^T) 
(p,\xi)}\circ d_{
(p,\xi)}\mathbb{LG}(\bn^T)\right)\\
&=  \det \left(\left(
 \begin{array}{cc}
 \langle -(\bn^T+\bn^S_{k-1})_{u_i},\bX_{u_j}\rangle & \langle -(\bn^T+\bn^S_{k-1})_{u_i},\bn^S_j \rangle \\
 ^{\atop{1 \leqslant i,j \leqslant s}} & ^{\atop{1 \leqslant i \leqslant s;1 \leqslant j \leqslant k-2}} \\
\bo_{(k-2) \times s}  &
-I_{(k-2)} \\
 \end{array}
 \right)
\left(
\begin{array}{cc}
g^{ij} & \bo \\
\bo & I_{(k-2)}
\end{array}
\right)\right)(u_0).
\end{split}
 \end{equation*}
Since $\ell_0\widetilde{\mathbb{LG}}(\bn^T,\bN^S)=\mathbb{LG}(\bn^T,\bN^S),$
we have
\begin{eqnarray*}
(\ell _0)_{u_i}\widetilde{\mathbb{LG}}(\bn^T,\bN^S)+\ell _0 \widetilde{\mathbb{LG}}(\bn^T,\bN^S)_{u_i}
&=& \mathbb{LG}(\bn^T,\bN^S)_{u_i}, \\
(\ell _0)_{\mu _j}\widetilde{\mathbb{LG}}(\bn^T,\bN^S)+\ell _0 \widetilde{\mathbb{LG}}(\bn^T,\bN^S)_{\mu_j}
&=& \mathbb{LG}(\bn^T,\bN^S)_{\mu_j}.
\end{eqnarray*}
Moreover, we have $\langle \widetilde{\mathbb{LG}}(\bn^T,\bN^S)(u_0,\mu_0),\bX_{u_i}(u_0)\rangle =
\langle \widetilde{\mathbb{LG}}(\bn^T,\bN^S)(u_0,\mu_0), \bn^S_j(u_0)\rangle =0.$
It follows that
\begin{equation*}\label{1012061439}
 \begin{split}
&\widetilde{K}_\ell (\bn^T)(p,\bxi)=
\det \left(-\Pi^t_{\widetilde{\mathbb{LG}}(n^T) 
(p,\xi)}\circ d_{
(p,\xi)}\widetilde{\mathbb{LG}}(\bn^T)\right)
\\
&=  \det \left(\left(
 \begin{array}{cc}
 \frac{1}{\ell_0}\langle -(\bn^T+\bn^S_{k-1})_{u_i},\bX_{u_j}\rangle & \frac{1}{\ell_0}\langle -(\bn^T+\bn^S_{k-1})_{u_i},\bn^S_j \rangle \\
 ^{\atop{1 \leqslant i,j \leqslant s}} & ^{\atop{1 \leqslant i \leqslant s;1 \leqslant j \leqslant k-2}} \\
\bo_{(k-2) \times s}  &
-\frac{1}{\ell_0}I_{(k-2)} \\
 \end{array}
 \right)
\left(
\begin{array}{cc}
g^{ij} & \bo \\
\bo & I_{(k-2)}
\end{array}
\right)\right)(u_0).
\end{split}
 \end{equation*}
\par
On the other hand, Corollary 3.2 implies that
\begin{equation*}\label{1012061440}
 \begin{split}
& K_{\ell}(\bn^T,\bxi)(p)  =
K_{\ell}(\bn^T,\bn^S_{k-1})(u_0) \\
& =\det ((\langle -(\bn^T+\bn^S_{k-1})_{u_i},\bX_{u_j})(g^{ij}))(u_0) \\
&=  \det \left(\left(
 \begin{array}{cc}
 \langle -(\bn^T+\bn^S_{k-1})_{u_i},\bX_{u_j}\rangle & \langle -(\bn^T+\bn^S_{k-1})_{u_i},\bn^S_j \rangle \\
 ^{\atop{1 \leqslant i,j \leqslant s}} & ^{\atop{1 \leqslant i \leqslant s;1 \leqslant j \leqslant k-2}} \\
\bo_{(k-2) \times s}  &
I_{(k-2)} \\
 \end{array}
 \right)
\left(
\begin{array}{cc}
g^{ij} & \bo \\
\bo & I_{(k-2)}
\end{array}
\right)\right)(u_0).
\end{split}
 \end{equation*}
Therefore we have the following theorem.
\begin{Th}
Under the same notations as those of the above paragraph, we have 
\begin{eqnarray*}
\widetilde{K}_\ell (\bn^T)(p_0,\bxi_0)&=&(-1)^{k-2}\mathcal{N}_h(\bn^T(u_0),\bxi_0)^{n-1}K_{\ell}(\bn^T,\bn^S)(u_0) \\
&=& (-\mathcal{N}_h(\bn^T(u_0),\bxi_0))^{k-2}\widetilde{K}_\ell (\bn^T,\bn^S)(u_0),
\end{eqnarray*}
where $p_0=\bX(u_0)$ and $\bn^S(u)$ is a local section of $N_1(M)[\bn^T]$ such that $\bn^S(u_0)=\bxi_0.$
\end{Th}
We have the following corollary of the above theorem.
\begin{Co}
The following conditions are equivalent{\rm :}
\par
{\rm (1)} $p_0=\bX(u_0)$ is a $(\bn^T_0,\bxi_0)$-parabolic point {\rm (} $K_\ell (\bn^T,\bn^S)(u_0)=0${\rm )}, 
\par
{\rm (2)} $\widetilde{K}_\ell(\bn^T)(p_0,\bxi_0)=0.$
\par\noindent
Here, $\bn^S(u)$ is a local section of $N_1(M)[\bn^T]$ such that $\bn^S(u_0)=\bxi_0.$

\end{Co}
\section{Lightcone height functions}
In order to investigate the geometric meanings of the normalized lightcone Lipschitz-Killing curvature of $N_1(M)[\bn^T]$, we  introduce
a family of functions on $M=\bX(U).$
\begin{Def}
{\rm
We define the family of {\it lightcone height functions} 
\[
H:U\times \mathbb{S}^{n-1}_+\lon \R
\] on $M=\bX(U)$ 
by
$
H(u,\bv)=\langle \bX(u),\bv\rangle.
$
We denote the Hessian matrix of the lightcone height function $h_{\sbv_0}(u)=H(u,\bv_0)$
at $u_0$ by ${\rm Hess}(h_{\sbv_0})(u_0).$
}\end{Def}
The following proposition characterizes the lightlike parabolic points and lightlike flat points
in terms of the family of lightcone height functions.
\begin{Pro}Let $H:U\times \mathbb{S}^{n-1}_+\lon \R$ be the family of lightcone height functions on M.
Then 
\par\noindent
{\rm (1)} $(\partial H/\partial u_i)(u_0,\bv_0)=0\ (i=1,\dots ,s)$ if and only if
there exists $\bxi_0\in N_1(M)_{p_0}[\bn^T]$ such that $\bv_0=\widetilde{\mathbb{LG}}(\bn^T)(p_0,\bxi_0),$ where
$p_0=\bX(u_0).$
\par
Suppose that $p_0=\bX(u_0)$, $\bv_0=\widetilde{\mathbb{LG}}(\bn^T)(p_0,\bxi_0)$. Then
\par\noindent
{\rm (2)} $p_0$ is a $(\bn^T_0,\bxi_0)$-parabolic point if and only if
${\rm det}\,{\rm Hess}(h_{\sbv_0})(u_0)=0,$ where $\bn^T_0=\bn^T(u_0),$
\par
\noindent
{\rm (3)} $p_0$ is a flat $(\bn^T_0,\bxi)$-umbilical point if and only if
${\rm rank}\,{\rm Hess}(h_{\sbv_0})(u_0)=0$,
\par\noindent
{\rm (4)} $u_0$ is a non-degenerate critical point of $h_{\sbv_0}$ if and only if $(p_0,\bxi _0)$ is a regular point of $\widetilde{\mathbb{LG}}(\bn^T).$
\end{Pro}
\demo
(1) Since $(\partial H/\partial u_i)(u_0,\bv_0)=\langle \bX_{u_i}(u_0),\bv_0),$
$(\partial H/\partial u_i)(u_0,\bv_0)=0\ (i=1,\dots ,s)$ if and only if $\bv_0\in N_{p_0}(M)$ and
$\bv_0\in \mathbb{S}^{n-1}_+.$
If $\langle \bv_0,\bn^T(u_0)\rangle =0,$ then $\bn^T(u_0)\in HP(\bv_0,0).$
But $HP(\bv_0,0)$ is a lightlike hyperplane. This fact contradicts to the fact that $\bn^T(u_0)$ is
timelike.
Thus, $\langle \bv_0,\bn^T(u_0)\rangle \not= 0.$
Then we can easily show that 
\[
\bxi_0=-\frac{1}{\langle\bn^T(u_0),\bv_0\rangle}\bv_0-\bn^T(u_0)\in N_1(M)_{p_0}[\bn^T].
\]
It follows that
\[
\bv_0=\widetilde{\bn^T(u_0)+\bxi_0}=\widetilde{\mathbb{LG}}(\bn^T) (p_0,\bxi_0).
\]
The converse also holds.
\par
For the proof of the assertions (2) and (3), as a consequence of Proposition 3.1, we have
\begin{eqnarray*}
{\rm Hess}(h_{\sbv_0})(u_0)&=&\left(\langle \bX_{u_iu_j}(u_0),\widetilde{\mathbb{LG}}(\bn^T)(p_0,\bxi_0)\rangle\right)
=\left( \frac{1}{\ell _0}\langle \bX_{u_iu_j}(u_0),\bn^T(u_0)+\bxi_0\rangle\right) \\
&=& \left(\frac{1}{\ell _0}\langle \bX_{u_i}(u_0),(\bn^T+\bn^S)_{u_j}(u_0)\rangle\right) \\
&=&\left(\frac{1}{\ell _0}\langle \bX_{u_i}(u_0),-\sum _{k=1}^s h^k_j(\bn^T,\bxi_0)(u_0)\bX_{u_k}(u_0)\rangle\right) \\
&=& \left(-\mathcal{N}_h(\bn^T(\underline{u_0}),\bxi _0)h_{ij}(\bn^T,\bxi_0)(u_0)\right),
\end{eqnarray*}
where $\bn^S(u)$ is a local section of $N_1(M)[\bn^T]$ such that $\bn^S(u_0)=\bxi_0.$
By definition,  $K_\ell (\bn^T,\bxi)(p_0)= 0$ if and only if ${\rm det}\, (h_{ij}(\bn^T,\bxi)(u_0))=0.$
The assertion (2) holds.
Moreover, $p_0$ is a flat $(\bn^T_0,\bxi_0)$-umbilical point if and only if $(h_{ij}(\bn^T,\bxi_0)(u_0))=O.$
So we have the assertion (3). 
\par
By the above calculation $u_0$ is a non-degenerate critical point of $h_{\sbv_0}$ if and only if
\[
\widetilde{K}_\ell (\bn^T,\bxi _0)(u_0)=\frac{\det\, \left(-\mathcal{N}_h(\bn^T(\underline{u_0}),\bxi _0)h_{ij}(\bn^T,\bxi_0)(u_0)\right)}{\det (g_{ij}(u_0))}\not=0.
\]
By Corollary 4.2, the last condition is equivalent to the condition $\widetilde{K}_\ell (\bn^T)(p_0,\bxi_0)\not= 0.$
By the definition of $\widetilde{K}_\ell (\bn^T)(p_0,\bxi_0)$, the above condition means that
$(p_0,\bxi _0)$ is a regular point of $\widetilde{\mathbb{LG}}(\bn^T).$
\enD

\section{The total absolute lightcone curvature}
\par
We have the following theorem.
\begin{Th}
Let $d\mathfrak{v}_{N_1(M) [\bn^T]}$ be the canonical volume form of $N_1(M)[\bn^T]$
and $d\mathfrak{v}_{\mathbb{S}^{n-1}_+}$ the canonical volume form of $\mathbb{S}^{n-1}_+.$
Then we have
\[
(\widetilde{\mathbb{LG}}(\bn^T )^*d\mathfrak{v}_{\mathbb{S}^{n-1}_+})_{( p,\xi ) }=|\widetilde{K}_\ell (\bn ^T)(p,\bxi)|
d\mathfrak{v}_{N_1(M) [\bn^T](p,\xi)}.
\]
\end{Th}
\demo
Without the loss of generality, we may assume that a point
$(p,\bxi)$ is a non-singular point of $\widetilde{\mathbb{LG}}(\bn^T,\bxi),$
We consider the same frame
$\{\bX_{u_1},\dots ,\bX_{u_s},\bn^T,\bn^S_1,\dots ,\bn^S_{k-1}\}$ as in the previous 
sections such that $\bn^S_{k-1}(u_0)=\bxi$ and $p=\bX(u_0).$
We also consider the local coordinate neighborhood 
$
U_{k-1}^+=\{(\mu_1,\dots ,\mu _{k-1})\in \mathbb{S}^{k-2}\ |\ \mu_{k-1}>0\ \},
$
of  $\mathbb{S}^{k-2}$,
so that we have $\mu_{k-1}=\sqrt{1-\sum_{j=1}^{k-2} \mu _j^2}$.
By the same calculations as just before Theorem 4.1, we have 
\begin{eqnarray*}
\mathbb{LG}( \bn^T,\bxi )_{u_i}(u_0)&=& (\bn^T+\bn^S_{k-1})_{u_i}(u_0) \\
&=&\langle \bn ^S_{u_i},\bn ^T\rangle(\bn^T-\bn^S_{k-1})(u_0)+\sum _{\ell =1}^{k-2}\langle (\bn^T+\bn^S_{k-1})_{u_i},\bn^S_\ell \rangle\bn^S_\ell(u_0) \\
&{}& -\sum_{j=1}^{s}
h_i^j(\bn^T,\bn^S )\bX _{u_j}(u_0).
\end{eqnarray*}
Therefore, we have
\begin{eqnarray*}
\langle\!\!\!\!\!\!\!\!\! &{}&\mathbb{LG}( \bn^T,\bxi )_{u_i}(u_0),\mathbb{LG}( \bn^T,\bxi )_{u_j}(u_0)
\rangle =
\langle \pi^{\tau}\circ (\bn^T+\bn^S_{k-1})_{u_i}(u_0),\pi^{\tau}\circ (\bn^T+\bn^S_{k-1})_{u_j}(u_0)\rangle \\
&{}& +\sum_{\ell =1}^{k-2} \langle (\bn^T+\bn^S_{k-1})_{u_i}(u_0),\bn^S_\ell(u_0)\rangle
\langle (\bn^T+\bn^S_{k-1})_{u_j}(u_0),\bn^S_\ell(u_0)\rangle.
\end{eqnarray*}
\par
It also follows from the calculations before Theorem 4.1 that
\begin{eqnarray*}
\langle \widetilde{\mathbb{LG}}( \bn^T,\bxi )_{u_i},\widetilde{\mathbb{LG}}( \bn^T,\bxi )_{u_j}\rangle
&=&\frac{1}{\ell _0^2}\langle \mathbb{LG}( \bn^T,\bxi )_{u_i},\mathbb{LG}( \bn^T,\bxi )_{u_j}\rangle,\\
\langle \widetilde{\mathbb{LG}}( \bn^T,\bxi )_{u_i},\widetilde{\mathbb{LG}}( \bn^T,\bxi )_{\mu_j}\rangle
&=&\frac{1}{\ell _0^2}\langle \mathbb{LG}( \bn^T,\bxi )_{u_i},\bn^S_j\rangle,\\
\langle \widetilde{\mathbb{LG}}( \bn^T,\bxi )_{\mu_i},\widetilde{\mathbb{LG}}( \bn^T,\bxi )_{\mu_j}\rangle
&=&\frac{1}{\ell _0^2}\langle \bn^S_i,\bn^S_j\rangle
\end{eqnarray*}
at $(u_0,\mu _0)\in U\times \mathbb{S}^{k-2}.$
We consider the matrix $A$ defined by
\begin{equation*}\label{1012061441}
 \begin{split}
A
&=  \left(
 \begin{array}{cc}
 \langle \widetilde{\mathbb{LG}}( \bn^T,\bxi )_{u_i},\widetilde{\mathbb{LG}}( \bn^T,\bxi )_{u_j}\rangle & \langle \widetilde{\mathbb{LG}}( \bn^T,\bxi )_{u_i},\widetilde{\mathbb{LG}}( \bn^T,\bxi )_{\mu_j}\rangle \\
 ^{\atop{1 \leqslant i,j \leqslant s}} & ^{\atop{1 \leqslant i \leqslant s;1 \leqslant j \leqslant k-2}} \\
\langle \widetilde{\mathbb{LG}}( \bn^T,\bxi )_{u_j},\widetilde{\mathbb{LG}}( \bn^T,\bxi )_{\mu_i}\rangle  &
\langle \widetilde{\mathbb{LG}}( \bn^T,\bxi )_{\mu_i},\widetilde{\mathbb{LG}}( \bn^T,\bxi )_{\mu_j}\rangle\\
^{\atop{1 \leqslant i \leqslant k-1;1\leqslant j \leqslant s}} & ^{\atop{1 \leqslant i,j \leqslant k-2}} \\
 \end{array}
 \right)(u_0).
\end{split}
 \end{equation*}
 By the previous calculation, we have
 \begin{equation*}\label{1012061442}
 \begin{split}
A
&=  \frac{1}{(\ell _0)^2}\left(
 \begin{array}{cc}
 \alpha _{ij} & \langle (\bn^T+\bn^S_{k-1})_{u_i},\bn^S_j \rangle \\
 ^{\atop{1 \leqslant i,j \leqslant s}} & ^{\atop{1 \leqslant i \leqslant s;1 \leqslant j \leqslant k-2}} \\
\langle (\bn^T+\bn^S_{k-1})_{u_j},\bn^S_i \rangle   &
\langle \bn^S_i,\bn^S_j \rangle\\
^{\atop{1 \leqslant i \leqslant k-1;1\leqslant j \leqslant s}} & ^{\atop{1 \leqslant i,j \leqslant k-2}} \\
 \end{array}
 \right)(u_0) \\
 &= \frac{1}{(\ell _0)^2}
 \left(
\begin{array}{cc}
 \alpha _{ij} & \langle (\bn^T+\bn^S_{k-1})_{u_i},\bn^S_j \rangle \\
 ^{\atop{1 \leqslant i,j \leqslant s}} & ^{\atop{1 \leqslant i \leqslant s;1 \leqslant j \leqslant k-2}} \\
\langle (\bn^T+\bn^S_{k-1})_{u_j},\bn^S_i \rangle   &
I_{k-1} \\
^{\atop{1 \leqslant i \leqslant k-1;1\leqslant j \leqslant s}} & ^{\atop{1 \leqslant i,j \leqslant k-2}}
\\
 \end{array}
 \right)(u_0),
\end{split}
 \end{equation*}
where 
\begin{eqnarray*}
\alpha_{ij}&=&
\langle \pi^{\tau}\circ (\bn^T+\bn^S_{k-1})_{u_i},\pi^{\tau}\circ (\bn^T+\bn^S_{k-1})_{u_j}\rangle \\
&{}& +\sum_{\ell =1}^{k-2} \langle (\bn^T+\bn^S_{k-1})_{u_i},\bn^S_\ell\rangle
\langle (\bn^T+\bn^S_{k-1})_{u_j},\bn^S_\ell\rangle .
\end{eqnarray*}
We consider a matrix
 \begin{equation*}\label{1012061443}
 \begin{split}
A_0 = \frac{1}{(\ell _0)^2}
 \left(
\begin{array}{cc}
\langle \pi^{\tau}\circ (\bn^T+\bn^S_{k-1})_{u_i},\pi^{\tau}\circ (\bn^T+\bn^S_{k-1})_{u_j}\rangle& \langle (\bn^T+\bn^S_{k-1})_{u_i},\bn^S_j \rangle \\
 ^{\atop{1 \leqslant i,j \leqslant s}} & ^{\atop{1 \leqslant i \leqslant s;1 \leqslant j \leqslant k-2}} \\
\bo_{(k-1)\times s}   &
I_{k-1} \\
\\
 \end{array}
 \right)(u_0).
\end{split}
 \end{equation*}
We denote that $A^j, A^j_0$ the $j$-the columns of the above two matrices.
Then we have the relation that
\[
A^j=A^j_0+\sum_{\ell=1}^{k-1}\langle (\bn^T+\bn^S_{k-1})_{u_j}, \bn^S_\ell\rangle A^{s+\ell}_0,
\]
for $j=1,\dots ,s.$ It follows that
\[
\det(A)= \det(A_0)=\frac{1}{(\ell_0)^{2(n-1)}}\det
\left(
\begin{array}{c}
\langle \pi^{\tau}\circ (\bn^T+\bn^S_{k-1})_{u_i},\pi^{\tau}\circ (\bn^T+\bn^S_{k-1})_{u_j}\rangle \\
 ^{\atop{1 \leqslant i,j \leqslant s}} 
 \end{array}
 \right)(u_0).
 \]
By Proposition 3.1, we have
$\pi^{\tau}\circ (\bn^T+\bn^S_{k-1})_{u_i}(u_0)=-\sum_{j=1}^s h^j_i(\bn^T,\bxi)(u_0)\bX_{u_j}(u_0),$
so that
\[
\langle \pi^{\tau}\circ (\bn^T+\bn^S_{k-1})_{u_i},\pi^{\tau}\circ (\bn^T+\bn^S_{k-1})_{u_j}\rangle(u_0)
=\sum _{\alpha,\beta} h^\alpha _i(\bn^T,\bxi)(u_0)h^\beta _j(\bn^T,\bxi)(u_0)g_{\alpha\beta}(u_0).
\]
It follows from Corollary 3.2 and Theorem 4.1 that
\[
\det(A)=\left(\frac{(-1)^{k-1}}{(\ell _0)^{(n-1)}}\right)^2(K_\ell(\bn^T,\bxi)(u_0))^2\det (g_{ij})
=(\widetilde{K}_\ell (\bn^T)(p,\bxi))^2\det(g_{ij}).
\]
This completes the proof.
\enD
\par
On the other hand, let $\overline{\bn}^T$ be another timelike unit normal future directed vector field along $M=\bX(U).$
Since the canonical action of $SO_0(1,n)$ on $\mathbb{H}^n(-1)$ is transitive, there exists 
$g\in SO_0(1,n)$ such that $g.\bn^T(u_0)=\overline{\bn}^T(u_0).$
Then we define a smooth mapping
\[
\Phi _g:N_1(M)_p [\bn^T]\lon N_1(M)_p [\overline{\bn}^T]
\]
by
$\Phi _g(p,\xi)=(p,g.\xi),$
where $p=\bX(u_0).$
By the definition of the canonical Riemannian metrics on $N_1(M)_p [\bn^T]$ and $N_1(M)_p [\overline{\bn}^T]$,
$\Phi _g$ is an isometry.
Therefore, we have
\[
\Phi^*_gd\mathfrak{v}_{N_1(M) [\overline{\bn}^T](p,g.\xi)}=d\mathfrak{v}_{N_1(M) [\bn^T](p,\xi)}.
\]
We define the $k-2$-dimensional lightcone unit sphere on the fibere as $\mathbb{S}^{k-2}_+(N(M)_p)=\mathbb{S}^{n-1}_+\cap N_p(M)$.
Then we have  $\widetilde{\mathbb{LG}}(\bn^T)(N_1(M)_p [\bn^T])\subset \mathbb{S}^{k-2}_+(N(M)_p).$
Moreover, we can easily show that 
\[
\widetilde{\mathbb{LG}}(\bn^T)|_{N_1(M)_p [\bn^T]}: N_1(M)_p [\bn^T]\lon \mathbb{S}^{k-2}_+(N(M)_p)
\]
is a diffeomorphism.
\par
There exists a differential form $d\sigma _{k-2}(\bn^T)$ of degree $k-2$ on
$N_1(M) [\bn^T]$ such that its restriction to a fiber is the volume element of 
the $k-2$-sphere.
We remark that 
\[
d\mathfrak{v}_{N_1(M) [\bn^T]}=d\mathfrak{v}_{M}\wedge d\sigma _{k-2}(\bn^T).
\]
Then we have the following key lemma:
\begin{Lem} Let $\bX:U\lon \R^{n+1}_1$ be a spacelike embedding with codimension $k$
and $\bn^T,\overline{\bn}^T$ be future directed unit timelike normal vector fields along
$M=\bX(U).$
For any $(p,\bxi)\in N_1(M)[\bn^T]$ with $p=\bX(u_0),$ $g\in SO_0(1,n)$ and $\Phi _g$ are given
in the previous paragraphs.
Then we have
\[
|\widetilde{K}_\ell (\bn^T)(p,\bxi)|d\sigma _{k-2}(\bn^T)_{\sbxi}
=
|\widetilde{K}_\ell (\overline{\bn}^T)(p,g.\bxi)|d\sigma _{k-2}(\overline{\bn}^T)_{g.\sbxi}
\]
and 
\[
\int _{N_1(M)_p [\bn^T]}|\widetilde{K}_\ell(\bn^T)(p,\bxi)| d\sigma _{k-2}(\bn^T)
=
\int _{N_1(M)_p [\overline{\bn}^T]}|\widetilde{K}_\ell(\overline{\bn}^T)(p,\overline{\bxi}) |d\sigma _{k-2}(\overline{\bn}^T).
\]
\end{Lem}
\demo
Under the previous notations, we have
\begin{eqnarray*}
&{}& \left(\widetilde{\mathbb{LG}}(\bn^T)|_{N_1(M)_p [\bn^T]}\right)^* d\frak{v}_{\mathbb{S}^{k-2}_+(N(M)_p)} \\
&{}& \hskip1cm=\left(\widetilde{\mathbb{LG}} (\bn^T )^* d\frak{v}_{\mathbb{S}^{n-1}_+}\right)|_{N_1(M)_p [\bn^T ]} 
=|\widetilde{K}_\ell (\bn^T)|d\sigma _{k-2}(\bn^T).
\end{eqnarray*}
\par
We remark that the canonical action of $SO_0(k-1)$ on $\mathbb{S}^{k-2}_+(N(M)_p)$
is transitive. For any $h\in SO_0(k-1)$, we denote that $\psi (h)(v)=h.v$ for $v\in \mathbb{S}^{k-2}_+(N(M)_p)$,
so that we have an isometry $\psi(h):\mathbb{S}^{k-2}_+(N(M)_p)\lon \mathbb{S}^{k-2}_+(N(M)_p).$
Thus, we have $\psi(h)^*d\frak{v}_{\mathbb{S}^{k-2}_+(N(M)_p)(v)}=d\frak{v}_{\mathbb{S}^{k-2}_+(N(M)_p)(h.v)}.$
\par
On the other hand, we have
\begin{eqnarray*}
&{}& \widetilde{\mathbb{LG}}(\overline{\bn}^T)|_{N_1(M)_p [\overline{\bn}^T]}\circ \Phi_g(p,\bxi)= \widetilde{\overline{\bn}^T(u)+g.\bxi}
\\
&{}& =\widetilde{g.(\bn^T (u )+\bxi )}=\psi (h)(\widetilde{(\bn^T (u )+\bxi )})=\psi (h)\circ \widetilde{\mathbb{LG}}(\bn^T)|_{N_1(M)_p [\bn^T]}(p,\bxi),
\end{eqnarray*}
for some $h\in SO_0(k-1).$
We set $\bv =\widetilde{(\bn^T (u )+\bxi )}=\widetilde{\mathbb{LG}}(\bn^T)|_{N_1(M)_p [\bn^T]}(p,\bxi)\in \mathbb{S}^{k-2}_+(N(M)_p).$
Then we have
\begin{eqnarray*}
&{}& (\widetilde{\mathbb{LG}}(\overline{\bn}^T)|_{N_1(M)_p [\overline{\bn}^T]}\circ \Phi_g)^*d\frak{v}_{\mathbb{S}^{k-2}_+(N(M)_p)(\sbv)} \\
&{}& = (\Phi _g)^*\left((\widetilde{\mathbb{LG}}(\overline{\bn}^T)|_{N_1(M)_p [\bn^T]})^*d\frak{v}_{\mathbb{S}^{k-2}_+(N(M)_p)(\sbv)}\right) \\
&{}&
=(\Phi _g)^*\left(|\widetilde{K}_\ell (\overline{\bn^T})|d\sigma _{k-2}(\overline{\bn}^T)_{\sbxi}\right)
\\
&{}&=|\widetilde{K}_\ell (\overline{\bn^T})(p,g.\bxi)|d\sigma _{k-2}(\overline{\bn}^T)_{g.\sbxi}.
\end{eqnarray*}
Moreover, we have
\begin{eqnarray*}
&{}& (\widetilde{\mathbb{LG}}(\bn^T)|_{N_1(M)_p [\bn^T]})^*\circ \psi(h)^*(d\frak{v}_{\mathbb{S}^{k-2}_+(N(M)_p)(\sbv)}) \\
&{}&=(\widetilde{\mathbb{LG}}(\bn^T)|_{N_1(M)_p [\bn^T]})^*d\frak{v}_{\mathbb{S}^{k-2}_+(N(M)_p)(h.\sbv)}
\\
&{}&=|\widetilde{K}_\ell (\bn^T)(p,\bxi)|d\sigma _{k-2}(\bn^T)_{\sbxi}.
\end{eqnarray*}
Since
$
\widetilde{\mathbb{LG}}(\overline{\bn}^T)|_{N_1(M)_p [\overline{\bn}^T]}\circ \Phi_g(p,\bxi)=\psi (h)\circ \widetilde{\mathbb{LG}}(\bn^T)|_{N_1(M)_p [\bn^T]}(p,\bxi),
$
we have
\[
|\widetilde{K}_\ell (\overline{\bn}^T)(p,g.\bxi)|d\sigma _{k-2}(\overline{\bn}^T)_{g.\sbxi}
=
|\widetilde{K}_\ell (\bn^T)(p,\bxi)|d\sigma _{k-2}(\bn^T)_{\sbxi}.
\]
Moreover, we have
\begin{eqnarray*}
\int _{N_1(M)_p [\bn^T]}|\widetilde{K}_\ell(\bn^T)(p,\bxi)| d\sigma _{k-2}(\bn^T)
&=&
\int _{\Phi _g(N_1(M)_p [\bn^T])}|\widetilde{K}_\ell(\bn^T)(p,g.\bxi)| d\sigma _{k-2}(g.\bn^T)\\
&=&
\int _{N_1(M)_p [\overline{\bn}^T]}|\widetilde{K}_\ell(\overline{\bn}^T)(p,\overline{\bxi}) |d\sigma _{k-2}(\overline{\bn}^T).
\end{eqnarray*}
This completes the proof.
\enD
We call the integral
\[
K_\ell^*(p)=\int _{N_1(M)_p [\bn^T]}|\widetilde{K}_\ell(\bn^T)(p,\bxi)| d\sigma _{k-2}(\bn^T)
\]
a {\it total absolute lightcone curvature} of $M$ at $p=\bX(u_0).$ 
In the global situation, we consider a closed orientable manifold $M$ with dimension $s$ and a spacelike immersion
$f:M\lon \R^{n+1}_1$.
We define the {\it total absolute lightcone curvature} of $M$ by the integral
\[
\tau _\ell (M,f)=\frac{1}{\gamma _{n-1}}\int _{M} K_\ell^*(p)d\mathfrak{v}_{M}=\frac{1}{\gamma _{n-1}}\int _{N_1(M)[\bn^T]}|\widetilde{K}_\ell(\bn^T)(p,\bxi)| d\mathfrak{v}_{N_1(M) [\bn^T]},
\]
where $\gamma _{n-1}$ is the volume of the unit $n-1$-sphere $\mathbb{S}^{n-1}.$
\section{The Chern-Lashof type theorem}
Let $f:M\lon \R^{n+1}_1$ be a spacelike immersion from an $s$-dimensional closed orientable manifold $M.$
We have the family of lightcone height functions $H:M\times \mathbb{S}^{n-1}_+\lon \R$ defined by
$H(x,\bv)=\langle f(x),\bv\rangle.$
By Proposition 5.1, $\bv\in \mathbb{S}^{n-1}_+$ is a critical value of $\widetilde{\mathbb{LG}}(\bn^T)$
if and only if there exists a point $p\in M$ such that $p$ is a degenerate critical point $h_{\sbv}$.
Therefore, we have the following proposition.
\begin{Pro}
The height function $h_{\sbv}$ is a Morse function if and only if $\bv$ is a regular value of
$\widetilde{\mathbb{LG}}(\bn^T)$.
\end{Pro}
\demo
By Proposition 5.1, $x\in M$ is a non-degenerate critical point of $h_{\sbv}$ if and only if
there exists $\bxi \in N_1(M)_{f(p)}[\bn^T]$ such that
$\bv=\widetilde{\mathbb{LG}}(\bn^T)(f(p),\bxi)$ and $(f(p),\bxi)$ is a regular point of
$\widetilde{\mathbb{LG}}(\bn^T).$
By definition, all critical points of a Morse function are non-degenerate, so that
the proof is completed.
\enD
Let $D\subset \mathbb{S}^{n-1}_+$ be the set of regular values of $\widetilde{\mathbb{LG}}(\bn^T)$.
Since $M$ is compact, $D$ is open and, by Sard's theorem, the complement of $D$ in $\mathbb{S}^{n-1}_+$ 
has null measure.
We define an integral valued function $\eta :D\lon \mathbb{N}$ by
\[
\eta (\bv)={\rm the\ number\ of\ elements\ of}\ \widetilde{\mathbb{LG}}(\bn^T)^{-1}(\bv),
\]
which turns out to be continuous.
\begin{Pro}
\[
\tau_\ell (M,f)
=\frac{1}{\gamma _{n-1}}\int _D \eta (\bv)d\mathfrak{v}_{\mathbb{S}^{n-1}_+}.
\]
\end{Pro}
\demo
For any $\bv\in D,$ there exists a neighborhood $U$ of $\bv$ in $D$ such that
$\widetilde{\mathbb{LG}}(\bn^T)^{-1}(U)$ is the disjoint union of connected open sets
$V_1,\dots ,V_k$, $k=\eta (\bv),$ on which $\widetilde{\mathbb{LG}}(\bn^T):V_i\lon U$ is a diffeomorphism.
By Theorem 4.1, we have
\[
\int _{V_i} |\widetilde{K}_\ell(\bn^T)| d\mathfrak{v}_{N_1(M) [\bn^T]}=
\int _{V_i} \widetilde{\mathbb{LG}}(\bn^T )^*d\mathfrak{v}_{\mathbb{S}^{n-1}_+}=
{\rm deg}\,(\widetilde{\mathbb{LG}}(\bn^T )|_{V_i})\int _{U} d\mathfrak{v}_{\mathbb{S}^{n-1}_+}
=\int _U d\mathfrak{v}_{\mathbb{S}^{n-1}_+}.
\]
Since $\gamma _{n-1}=\int _{\mathbb{S}^{n-1}_+} d\mathfrak{v}_{\mathbb{S}^{n-1}_+}$ and $\widetilde{K}_\ell(\bn^T)$ is zero
at a singular point of $\widetilde{\mathbb{LG}}(\bn^T ),$ we have
\[
\frac{1}{\gamma _{n-1}}\int _{N_1(M)[\bn^T]}|\widetilde{K}_\ell(\bn^T)| d\mathfrak{v}_{N_1(M) [\bn^T]}
=
\frac{1}{\gamma _{n-1}}\int _D \eta (\bv)d\mathfrak{v}_{\mathbb{S}^{n-1}_+}.
\]
\enD
\par
We recall that the Morse number of a compact manifold $M$, $\gamma (M)$, is defined
to be the minimum number of critical points for any Morse function $\phi :M\lon \R.$
\begin{Th}[The Chern-Lashof type theorem]
Let $f:M\lon \R^{n+1}_1$ be a spacelike immersion of a compact $s$-dimensional manifold $M.$
Then
\par
{\rm (1)} $\tau_\ell (M,f)\geq \gamma (M)\geq 2,$
\par
{\rm (2)} If $\tau _\ell (M,f)< 3,$ then $M^s$ is homeomorphic to an $s$-sphere.
\end{Th}
\demo
Since each Morse function $h_{\sbv}$ certainly satisfies $\eta (\bv)\geq \gamma (M),$
we have $\tau _\ell (M,f)\geq \gamma (M).$
Since $M$ is compact, there exist at least two critical points for any smooth function
on $M$, so that $\gamma (M)\geq 2.$
If  $\tau _\ell (M,f)< 3,$ there must be a set $U$ of positive measure on which $\eta (\bv)=2.$
So there is a non-degenerate $h_{\sbv}$ with two critical points, and $M$ is homeomorphic to
$\mathbb{S}^s$ by Reeb's theorem (see, \cite{M63}).
\enD
\begin{Def}{\rm
If $\tau _\ell (M,f)=\gamma (M)$, then every non-degenerate lightlike height function $h_{\sbv}$ has the minimum number of critical points allowed by the Morse inequalities.
In this case we say that $f$ is a {\it lightlike-tight spacelike immersion} (or, simply,
{\it L-tight spacelike immersion}).
}
\end{Def}
In \S 9, we consider the problem to characterize the L-tightness for
spacelike immersed spheres.

\section{Codimension two spacelike submanifolds}
In the case when $s=n-1,$ 
$N_1(M)[\bn^T]$ is a double covering of $M.$ 
If $M$ is orientable, we can choose global section $\bsigma (p)=(p,\bn^S(p))$ of $N_1(M)[\bn^T].$
Let $\pi :\R^{n+1}\lon \R^n_0$ be the canonical projection
defined by $\pi (x_0,x_1,\dots ,x_n)=(0,x_1,\dots ,x_n),$
where $\R^n_0$ is the Euclidean space given by $x_0=0.$
Since ${\rm Ker}\, d\pi _{f(p)}$ is a timelike one-dimensional subspace of $\R^{n+1}_1$
and $\bn^S$ is spacelike, $d\pi _{f(p)}(\bn^S(p))$ is transverse to $\pi\circ f(M)$ at
$p\in M.$ Therefore, if $M$ is closed and $f:M\lon \R^{n+1}_1$ is a spacelike embedding such that
$\pi\circ f:M\lon \R^n_0$ is an embedding, then we can choose the direction of $\bn^S$ such that
$d\pi\circ \bn^S$ points the direction to the outward of $\pi\circ f(M).$

In \cite{IzCar07} it has been shown that $\widetilde{(\bn^T(p)\pm\bn^S(p))}$ is independent of the choice of
$\bn^T.$ 
Therefore, we have the global lightcone Gauss map 
\[
\widetilde{\mathbb{LG}}_{\pm} : M\lon \mathbb{S}^{n-1}_+
\]
defined by $\widetilde{\mathbb{LG}}_{\pm}(p)=\widetilde{(\bn^T(p)\pm\bn^S(p))}.$
Moreover, we have defined the normalized lightcone Lipschitz-Killing (i.e., Gauss-Kronecker) curvature $\widetilde{K}_\ell^{\pm}(p)=\widetilde{K}_\ell(\bn^T,\pm\bn^S)(p)$ of $M$ in \cite{IzCar07} .
Since $\widetilde{\mathbb{LG}}(\bn^T)(p,\pm\bn^S(p))=\widetilde{\mathbb{LG}}_{\pm}(p),$
we have
\[
\widetilde{K}_\ell^{\pm}(p)=\widetilde{K}_\ell(\bn^T,\pm\bn^S)(p)=\widetilde{K}_\ell(\bn^T)(p,\pm\bn^S(p)).
\]
In \cite{IzCar07} we have shown the following Gauss-Bonnet type theorem:
\begin{Th}[\cite{IzCar07}]
Suppose that $M$ is a closed orientable $n-1$-dimensional manifold, $n-1$ is even and $f:M\lon\R^{n+1}_1$ is
a spacelike embedding. Then 
\[
\int_{M}\widetilde{K}^\pm _\ell d\mathfrak{v}_{M^{n-1}}=\frac{1}{2}\gamma _{n-1}\chi (M),
\]
where $\chi (M^{n-1})$ is the Euler characteristic of $M^{n-1}$.
\end{Th}
In order to prove the above theorem, it has been shown in \cite{IzCar07} that
$\widetilde{K}^\pm _\ell d\mathfrak{v}_{M}=(\widetilde{\mathbb{LG}}_{\pm})^*d\mathfrak{v}_{\mathbb{S}^{n-1}_+}.$
Let $D^\pm\subset \mathbb{S}^{n-1}_+$ denote the set of regular value of $\widetilde{\mathbb{LG}}_{\pm}$
and $D=D^+\cap D^-.$
We define a mapping $\eta ^{\pm}:D\lon \mathbb{N}$
by
\[
\eta ^\pm (\bv)={\rm the\ number\ of\ elements\ of}\ (\widetilde{\mathbb{LG}}_{\pm})^{-1}(\bv).
\]
We have the following proposition:
\begin{Pro}
Suppose that $M$ is a closed orientable $n-1$-dimensional manifold and $f:M\lon\R^{n+1}_1$ is
a spacelike embedding. Then 
\[
\int_{M}|\widetilde{K}^\pm _\ell |d\mathfrak{v}_{M}=\int _D \eta ^{\pm}(\bv)d\mathfrak{v}_{\mathbb{S}^{n-1}_+}.
\]
\end{Pro}
\demo
Since $\widetilde{K}^\pm _\ell d\mathfrak{v}_{M}=(\widetilde{\mathbb{LG}}_{\pm})^*d\mathfrak{v}_{\mathbb{S}^{n-1}_+},$
we can prove by exactly the same arguments as those in the proof of Proposition 7.2.
\enD
\begin{Th}
Suppose that $M$ is a closed orientable $n-1$-dimensional manifold and $f:M\lon\R^{n+1}_1$ is
a spacelike embedding  such that $\pi\circ f$ is an embedding. Then 
\[
\int_{M}|\widetilde{K}^\pm _\ell |d\mathfrak{v}_{M}\geq \gamma_{n-1}.
\]
The equality holds if and only if $\widetilde{\mathbb{LG}}_{\pm}$ is bijective on the regular values.
\end{Th}
\demo
Since $\pi\circ f$ is an embedding, we can choose the vector field $\bn^S$ along $M$ such that
$d\pi\circ \bn^S$ is a transversal inward vector filed over $\pi\circ f(M)$ in $\R^n_0.$
It is enough to show that both of $\widetilde{\mathbb{LG}}_{\pm}$ are surjective onto $D.$
By Proposition 3.3, $p\in M$ is a critical point of the lightcone height function $h_{\sbv}$ if and only if
$\bv =\widetilde{\mathbb{LG}}(\bn^T)(p,\bxi)$ for some $\bxi \in N_1(M)_p[\bn^T].$
Since the codimension of $M$ is two, the last condition is equivalent to the condition $\bv =\widetilde{\mathbb{LG}}(\bn^T)(p,\bxi)=\widetilde{\mathbb{LG}}_+ (p)$
or $\bv =\widetilde{\mathbb{LG}}(\bn^T)(p,-\bxi)=\widetilde{\mathbb{LG}}_- (p).$
For any $\bv\in \mathbb{S}^{n-1}_+,$ there exists the maximum point $p_0$ and the minimum point $q_0$ of the
lightcone height function $h_{\sbv}$ on the compact manifold $M.$
These points are critical points of $h_{\sbv},$ so that $\bv=\widetilde{\mathbb{LG}}_+(p_0)$ or $\bv=\widetilde{\mathbb{LG}}_-(p_0)$ ( and $\bv=\widetilde{\mathbb{LG}}_+(q_0)$ or $\bv=\widetilde{\mathbb{LG}}_-(q_0)$). It is enough to show that $\widetilde{\mathbb{LG}}_+(p_0)\not= \widetilde{\mathbb{LG}}_+(q_0).$
Suppose that $\widetilde{\mathbb{LG}}_+(p_0)= \widetilde{\mathbb{LG}}_+(q_0).$
We define a function $\widetilde{h}_{\sbv} :\R^4_1\lon \R$ by $\widetilde{h}_{\sbv}(\bx)=\langle \bv,\bx\rangle.$
It follows that $\widetilde{h}_{\sbv}\circ f(p)=h_{\sbv}(p).$ We distinguish two cases.
\par
(i) If $\bv=\widetilde{\mathbb{LG}}_+(p_0),$ then we have $\bv=\widetilde{\mathbb{LG}}_+(q_0).$
We consider the line from $f(q_0)$ directed by $-\bn^S(q_0),$ parametrized by
\[
\bgamma_{q_0}(t)=f(q_0)-t\bn^S(q_0).
\]
Then we have
\begin{eqnarray*}
\frac{d\widetilde{h}_{\sbv}\circ\bgamma _{q_0}}{dt}(t) &=&
\langle -\bn^S(q_0),\bv\rangle =\langle -\bn^S(q_0),\widetilde{\mathbb{LG}}_+(q_0)\rangle \\
&=& \left\langle -\bn^S(q_0),\frac{1}{\ell _0^+(q_0)}(\bn^T(q_0)+\bn^S(q_0)\right\rangle 
= -\frac{1}{\ell _0^+(q_0)}<0.
\end{eqnarray*}
It follows that $\widetilde{h}_{\sbv}\circ\bgamma _{q_0}(t)$ is strictly decreasing.
Since $q_0$ is the minimum point of $h_{\sbv}$ and $f(q_0)=\bgamma _{q_0}(0),$
$\bgamma _{q_0}(t)\notin f(M)$ for any $t>0.$
Thus, we have $\pi\circ \bgamma _{q_0}(t)\notin \pi\circ f(M)$ for any $t>0.$
Since $\pi\circ \bgamma _{q_0}$ is a line in $\R^n_0,$ there exists a positive real number $\tau $
such that $\pi\circ \bgamma _{q_0}(\tau )$ is in the outside of $\pi\circ f(M).$
On the other hand, since $d\pi\circ \bn^S$ is an inward transversal vector field along $\pi\circ f(M)$ in
$\R^n_0,$ there exists a sufficiently small $\varepsilon >0$ such that
$\pi\circ \bgamma _{q_0}(\varepsilon )$ is in the inside of $\pi\circ f(M).$ By the Jordan-Brouwer separation theorem, there exists a real number $t_0>0$ such that $\pi\circ \bgamma _{q_0}(t_0 )\in \pi\circ f(M).$ 
This is a contradiction.
\par
(ii) If $\bv=\widetilde{\mathbb{LG}}_-(p_0),$ then we also consider the line from $f(p_0)$ defined by
\[
\bgamma_{p_0}(t)=f(p_0)-t\bn^S(p_0).
\]
Then we have
\begin{eqnarray*}
\frac{d\widetilde{h}_{\sbv}\circ\bgamma _{p_0}}{dt}(t) &=&
\langle -\bn^S(p_0),\bv\rangle =\langle -\bn^S(p_0),\widetilde{\mathbb{LG}}_-(p_0)\rangle \\
&=& \left\langle -\bn^S(p_0),\frac{1}{\ell _0^+(p_0)}(\bn^T(p_0)-\bn^S(p_0)\right\rangle 
= \frac{1}{\ell _0^+(q_0)}>0,
\end{eqnarray*}
so that $\widetilde{h}_{\sbv}\circ\bgamma _{p_0}(t)$ is strictly increasing.
Since $p_0$ is the maximum point of $h_{\sbv}$ and $f(p_0)=\bgamma _{p_0}(0),$
$\bgamma _{p_0}(t)\notin f(M)$ for any $t>0.$
By exactly the same reason as in the case (i), there exists a real number $t_0>0$ such that
$\bgamma _{p_0}(t_0)\in \pi\circ f(M)$.
This is a contradiction.
\enD
\begin{Def}{\rm
We define the {\it total absolute lightcone curvature} of a spacelike embedding $f:M\lon\R^{n+1}_1$
from a closed orientable $n-1$-dimensional manifold by
\[
\tau ^\pm_\ell (M,f)=\frac{1}{\gamma _{n-1}}\int_{M}|\widetilde{K}^\pm _\ell |d\mathfrak{v}_{M}.
\]
}
\end{Def}
\par
We remark that we have the following weaker inequality from Theorem 7.3:
\[
\tau ^+_\ell (M,f)+\tau ^-_\ell (M,f)=\tau _\ell (M,f)\geq 2.
\]
There are examples such that
\[
\tau ^+_\ell (M,f)\not=\tau ^-_\ell (M,f)
\]
(see Subsection 10.2).
\par
For an even dimensional manifold $M,$ we have the following theorem.
\begin{Th} Let $f:M\lon \R^{n+1}_1$ be a spacelike embedding from a
closed orientable $n-1$-dimensional manifold.
Suppose $n$ is an odd number.
Then we have
\[
\int_{M}|\widetilde{K}^\pm _\ell |d\mathfrak{v}_{M}\geq \frac{1}{2}\gamma _{n-1}(4-\chi (M)),
\]
\end{Th}
where $\chi (M)$ is the Euler characteristic of $M.$
\demo
In order to avoid the confusion, we only give a proof for $\widetilde{K}^+ _\ell.$ 
Consider the lightcone Gauss map
$\widetilde{\mathbb{LG}}_+:M\lon \mathbb{S}^{n-1}_+.$
We define $M^+=\{p\in M \ |\ \widetilde{K}_\ell^+ >0\}$ and
$M^-=\{p\in M \ |\ \widetilde{K}_\ell^+ <0\}.$
Then we can write
\[
\int_{M}|\widetilde{K}^+ _\ell |d\mathfrak{v}_{M}=\int_{M^+}\widetilde{K}^+_\ell d\mathfrak{v}_{M} - \int_{M^-}\widetilde{K}^+ _\ell d\mathfrak{v}_{M}
\]
and
\[
\int_{M}\widetilde{K}^+_\ell d\mathfrak{v}_{M}=
\int_{M^+}\widetilde{K}^+_\ell d\mathfrak{v}_{M} + \int_{M^-}\widetilde{K}^+ _\ell d\mathfrak{v}_{M}.
\]
By Theorem 8.1 and the above equations, we have
\[
\int_{M}|\widetilde{K}^+ _\ell |d\mathfrak{v}_{M}=
2\int_{M^+}\widetilde{K}^+  _\ell d\mathfrak{v}_{M}-\frac{1}{2}\gamma _{n-1}\chi (M).
\]
Thus, it is enough to show that
\[
\int_{M^+}\widetilde{K}^+ _\ell d\mathfrak{v}_{M}\geq \gamma _{n-1}.
\]
Let $M_0,M_1,M_2,M^+_2$ be the subsets of $M$ defined by
$M_0=(\widetilde{K}^+ _\ell)^{-1}(0),$ $M_1=\{p\in M\setminus M_0\ |\ \exists q\in M_0\ with \ \widetilde{\mathbb{LG}}_+ (q)=\widetilde{\mathbb{LG}}_+ (p)\ \},$
$M_2=M\setminus (M_0\cup M_1)$ and $M^+_2=M^+\cap  M_2.$
Since $M_0$ is the singular set of $\widetilde{\mathbb{LG}}_+,$ $\widetilde{\mathbb{LG}}_+(M_0)$ hs measure zero by Sard's Theorem
and also $\widetilde{\mathbb{LG}}_+(M_0)\cup \widetilde{\mathbb{LG}}_+(M_1)$ is a mesure
zero set in $S^2_+.$ For any $\bv\in S^2_+\setminus (\widetilde{\mathbb{LG}}_+(M_0)\cup \widetilde{\mathbb{LG}}_+(M_1)),$ the lightcone height function $h_{\sbv}$ has at least two critical points: a maximum and a minimum.
In \cite{IzCar07}, it was shown that
\[
\widetilde{K}^+ _\ell (p)=\frac{{\rm det}\, {\rm Hess}\,(h_{\sbv}(p))}{{\rm det}\, (g_{ij}(p))},
\]
where $\bv =\widetilde{\mathbb{LG}}_+(p).$
Since $\bv$ is a regular value of $\widetilde{\mathbb{LG}}_+,$
$h_{\bv}$ has a Morse-type singular point with index $0$ or $n-1$ at the minimum point and the maximum point.
The lightcone Gauss-Kronecker curvature $\widetilde{K}^+ _\ell $ is positive at such points, so that $\widetilde{\mathbb{LG}}_+ |_{M^+}$ is surjective.
For the case of $\widetilde{K}^- _\ell $, we can show the assertion by exactly the same arguments as the above case.
\enD
As a special case for $n=3,$ we have the following corollary.
\begin{Co}
For a spacelike embedding $f:M\lon \R^4_1$ from a closed orientable surface $M,$
we have
\[
\int_{M}|\widetilde{K}^\pm _\ell |d\mathfrak{v}_{M}\geq 2\pi (4-\chi (M)),
\]
\end{Co}
We define the {\it lightcone mean curvature} of $M$ at $p$ by 
\[
\widetilde{H}^\pm _\ell(p)=\frac{1}{2}{\rm Trace}\,\widetilde{S}^\pm_p=\frac{1}{2}(\widetilde{\kappa}^\pm_1(p)+\widetilde{\kappa}^\pm_2(p)),
\]
where $\widetilde{S}^\pm_p=\widetilde{S}(\bn^T,\pm\bn^S)_p.$
Then we have the following proposition.
\begin{Pro}
For a spacelike embedding $f:M\lon \R^4_1$ from a closed orientable surface $M,$
we have
\[
\int_{M} (\widetilde{H}^\pm _\ell)^2 d\mathfrak{v}_{M}\geq 4\pi.
\]
The equality holds if and only if $f:M\lon \R^4_1$ is totally umbilical with a non-zero normalized principal curvature.
\end{Pro}
\demo
Since $\widetilde{H}^\pm _\ell(p)=(\widetilde{\kappa}^\pm_1(p)+\widetilde{\kappa}^\pm_2(p))/2$
and $\widetilde{K}^\pm _\ell=\widetilde{\kappa}^\pm_1(p)\widetilde{\kappa}^\pm_2(p),$
we have
\[
(\widetilde{H}^\pm _\ell)^2 (p)-\widetilde{K}^\pm _\ell(p)=
\frac{1}{4}(\widetilde{\kappa}^\pm_1(p)-\widetilde{\kappa}^\pm_2(p))^2\geq 0.
\]
It follows that 
\[
\int_{M} (\widetilde{H}^\pm _\ell)^2 d\mathfrak{v}_{M}\geq \int_{M^+} (\widetilde{H}^\pm _\ell)^2 d\mathfrak{v}_{M}\geq \int_{M^+} \widetilde{K}^\pm _\ell d\mathfrak{v}_{M}.
\]
By the assertion in the proof of Theorem 8.4, we have
\[
 \int_{M^+} \widetilde{K}^\pm _\ell d\mathfrak{v}_{M}\geq \gamma _2=4\pi.
 \]
 The equality holds if and only if 
 \[
 \int_{M^+} \left((\widetilde{H}^\pm _\ell )^2 -\widetilde{K}^\pm _\ell\right) d\mathfrak{v}_{M}=0.
 \]
 This means that $\widetilde{\kappa}^\pm_1(p)=\widetilde{\kappa}^\pm_2(p)$ for any $p\in M,$
 so that $M$ is totally lightcone umbilical. This completes the proof.
\enD
\par
\begin{Rem}\rm
(1) In \cite{Izu-Pei-Romero2} it was shown that there exists a parallel timelike future directed
unit normal vector field $\bn^T$ along $f:M\lon \R^4_1$ and totally umbilical with a non-zero lightcone principal curvature  if and  only if  $M$ is embedded in the lightcone.
 It is well known that if a compact surface $M$ is embedded in the lightcone, it is homeomorphic to
a sphere.
In this case the normalized lightcone principal curvature is constant, but
the lightcone principal curvature is not constant.
So, the surface $f(M)$ is not necessarily a round sphere.
\par
On the other hand,
suppose that  $f(M)$ is in the Euclidean space or the hyperbolic space.
Since the intersection of the lightcone with Euclidean space or the hyperbolic space is
a round sphere, the equality of the above theorem holds
if and only if $f(M)$ is a round sphere.
\par\noindent
(2) In the first draft of this paper, we proposed the lightcone version of the Willmore conjecture.
However, the anonymous referee has pointed out there exists a spacelike immersion $f:T\lon \R^4_1$
from the torus such that
\[
\int_{T} (\widetilde{H}^\pm _\ell)^2 d\mathfrak{v}_{T}< 2\pi^2.
\]
If $T$ is immersed into the Euclidean space $\R^3_0,$
then we have the original Willmore conjecture (cf.\S 10).
Recently, the Willmore conjecture has been proved by F. C. Marques and A. Neves in \cite{M-N}.
Moreover, if $T$ is immersed into the hyperbolic space $\mathbb{H}^3(-1),$ we have
the horospherical Willmore conjecture (cf., \S 10).
Therefore we have the following new problem.
\vskip3pt
\noindent
{\bf Problem.} What value is the lower bound of the lightcone Willmore energies for spacelike tori in $\R^4_1$? 
\vskip3pt
\end{Rem}
 
\section{Lightlike tight spacelike spheres}
In this section we consider the characterizations of L-tightness for spacelike spheres. 
Let $f:M\lon \R^{n+1}_1$ be a spacelike immersion of a closed orientable manifold $M.$
We remind the reader that $f$ is called an L-tight if every non-degenerate lightcone
height function $h_{\sbv}$ has the minimum number of critical points required by the Morse
inequalities. If $M$ is homeomorphic to a sphere, then the Morse number $\gamma (M)$ is equal to 2.
We have the following theorem.
\begin{Th} Let $f:M\lon \R^{n+1}_1$ be a spacelike immersion of a closed orientable
manifold $M$. Then the following conditions are equivalent:
\par
{\rm (1)} $M$ is homeomorphic to a sphere and $f$ is L-tight,
\par
{\rm (2)} $\tau _\ell (M,f)=2.$
\end{Th}
\demo
We use the function $\eta :D\lon \mathbb{N}$ defined before Proposition 7.2 in \S 7.
Here, $D$ is the regular value set of $\widetilde{\mathbb{LG}}(\bn^T).$
Since $M$ is compact, $D$ is open and $\mathbb{S}^{n-1}_+\setminus D$ has null measure by
the Sard theorem.
By Proposition 7.2, $\tau _\ell (M,f)=2$ if and only if $\eta (\bv)=2.$
This condition is equivalent to the following condition:
\par
$(*)$ The lightcone Gauss map of $N_1(M)[\bn^T]$ takes every regular value exactly twice.
\par\noindent
Suppose that the condition (1) holds. Then $\gamma (M)=2$.
Since $f$ is L-tight, the lightcone height function $h_{\sbv}$ for $\bv\in D$ has exactly 
$\gamma (M)=2$ non-degenerate critical points.
This is equivalent to the condition $(*)$.
For the converse, suppose that the condition $(*)$ holds.
Then $h_{\sbv}$ for $\bv\in D$ has exactly $2$ non-degenerate critical points,
so that $f$ is L-tight.
By the assertion (2) of Theorem 7.3, $M$ is homeomorphic to a sphere.
This completes the proof.
\enD
\par
By the above theorem, if $M$ is a sphere, $\tau _\ell (\mathbb{S}^s,f)=2$ if and only if $f$ is L-tight.
In order to give a further characterization, we introduce the following notion:
Let $V$ be a codimension two spacelike affine subspace of $\R^{n+1}_1.$
We define
$\overline{V}$ as a spacelike subspace parallel to $V$.
Since $\overline{V}^\perp$ is a Lorentz plane, there exists a pseudo-orthonormal basis
$\{\bv^T,\bv^S\}$ of $\overline{V}^\perp$ then we have lightlike vectors $\bv^+=\bv^T+\bv^S, \bv^-=\bv^T-\bv^S.$
There exists $p\in \R^{n+1}_1$ such that $V=p+\overline{V}.$
For any $\bw\in \overline{V},$ $\langle p+\bw,\bv^\pm\rangle=\langle p,\bv^\pm\rangle=c^\pm$ are
constant numbers.
We consider lightlike hyperplanes $HP(\bv^\pm,c^\pm).$
Then we have 
\[
V=HP(\bv^+,c^+)\cap HP(\bv^-,c^-).
\]
For a point $p\in M,$ we say that 
a codimension two spacelike affine subspace $V$ is a {\it codimension two spacelike tangent space}
if $T_pM\subset \overline{V}.$
Moreover, each one of $HP(\bv^\pm,c^\pm)$ is said to be a {\it tangent lightlike hyperplane\/} of 
$M$ at $p.$
Let $K$ be a subset of $\R^{n+1}_1$.
A hyperplane $HP$ through a point $\bx\in K$ is called a {\it support plane\/} of $K$ if $K$ lies entirely in
one of the closed half-spaces determined by $HP.$ The half-space is called a {\it support half-space\/}. 
Let $M$ be a compact orientable $n-1$-dimensional manifold.
Then we have unique two lightlike tangent hyperplanes of $f(M)$ at each point $p\in M.$
These hyperplanes are 
$HP(\bv^\pm,c^\pm)$, where $\bv^\pm=\bn^T(p)\pm\bn^S(p)$ and $c^\pm=\langle f(p),\bn^T(p)\pm\bn^S(p)\rangle.$
In this case, we say that $f(M)$ is {\it lightlike convex\/} (or, {\it L-convex} in short) if
for each $p\in M,$ the lightlike tangent hyperplanes of $f(M)$ at $f(p)$ are support planes of $f(M).$
\par
We consider the case that $M$ is a sphere.
Let $f:\mathbb{S}^s\lon \R^{n+1}_1$ be a spacelike immersion.
If $s=n-1$, we have the following theorem.
\begin{Th}
Let $f:\mathbb{S}^{n-1}\lon \R^{n+1}_1$ be a spacelike embedding.
Then the following conditions are equivalent:
\par
{\rm (1)} $f$ is L-convex,
\par
{\rm (2)} $\tau _\ell (\mathbb{S}^{n-1},f)=2$,
\par
{\rm (3)} $f$ is L-tight.
\par\noindent
Generally the following condition {\rm (4)} implies the condition {\rm (2)}.
If we assume that $n$ is odd or $\pi\circ f:M\lon \R^n_0$ is an embedding,
then the condition {\rm (2)} implies the condition {\rm (4)}.
\par
{\rm (4)} $\tau ^+_\ell (\mathbb{S}^{n-1},f)=\tau ^-_\ell (\mathbb{S}^{n-1},f)=1.$
\end{Th}
\demo
By Theorem 9.1, the conditions (2) and (3) are equivalent.
By Theorem 8.4, the condition (2) implies (4) for the case when $n$ is odd.
If $\pi\circ f$ is an embedding, Theorem 8.3 asserts that the condition (2) implies (4) 
even for the case when $n$ is even.
It is trivial that the condition (4) implies the condition (2).
\par
We now give a proof that the conditions (1) and (3) are equivalent.
Suppose that $f$ is L-tight. If $f$ is not L-convex,
then there exists $p\in \mathbb{S}^{n-1}$ and $\bv\in \mathbb{S}^{n-1}_+$ such that one of the tangent
lightlike hyperplanes at $p$ separates $f(\mathbb{S}^{n-1})$ into two parts.
Then we have $\bv=\widetilde{\mathbb{L}}^+ (p)$ or $\bv=\widetilde{\mathbb{L}}^- (p)$ (i.e., $p$ is a
critical point of $h_{\sbv}$).
If $p$ is a non-degenerate critical point, it contradicts to the assumption that $f$ is
L-tight.
If $p$ is a degenerate critical point, under a small perturbation of $\bv\in \mathbb{S}^{n-1}_+,$
we have a non-degenerate critical point which is neither the maximum nor the minimum point.
This also contradicts to the assumption that $f$ is L-tight.
We now suppose that f is not L-tight.
Then there exists a non-degenerate lightcone height function $h_{\sbv}$ which as at least
three critical points. If necessary, under a small perturbation of $\bv\in \mathbb{S}^{n-1}_+,$ all critical
values of $h_{\sbv}$ are different. It follows that there exists a critical point $p\in \mathbb{S}^{n-1}$ of
$h_{\sbv}$ such that neither the maximum nor the minimum point of $h_{\sbv}.$
This means that one of the tangent lightlike hyperplanes of $f(\mathbb{S}^{n-1})$ locally separates $f(\mathbb{S}^{n-1})$
into at least two parts.
Therefore, $f$ is not lightlike convex.
\enD
\par
We now consider the case when $n+1-s>2.$ 
For any $(p,\xi)\in N_1(M)[\bn^T],$ we consider the lightlike tangent hyperplanes
$HP(\bv^\pm _p,c^\pm),$ where $\bv^\pm_p=\bn^T(p)\pm\xi$ and $c^\pm=\langle f(p),\bv^\pm_p\rangle.$
We denote that $T_SM[\bn^T,\xi]_p=HP(\bv^+_p,c^+)\cap HP(\bv^-_p,c^-),$
which is called a {\it spacelike tangent affine space with codimension two\/} of $f(M)$ at $p\in M.$
We also define 
\begin{eqnarray*}
F_\ell (\bn^T(p),\pm\xi)&=&\{\bx\in \R^{n+1}_1\ |\ \langle \bx-f(p),\bv^\pm_p\rangle\leq 0\ \} \\
P_\ell (\bn^T(p),\pm\xi)&=&\{\bx\in \R^{n+1}_1\ |\ \langle \bx-f(p),\bv^\pm_p\rangle\geq 0\ \}.
\end{eqnarray*}
We call $F_\ell (\bn^T(p),\xi)$ (respectively, $P_\ell (\bn^T(p),\pm\xi)$) the {\it future regions\/} (respectively,
the {\it past regions\/}) with respect to $(\bn^T(p),\pm\xi).$
We have a closed subset
\[
S(\bn^T(p),\xi)=
\R^{n+1}_1\setminus {\rm Int}\left((F_\ell (\bn^T(p),+\xi)\cap F_\ell(\bn^T(p),-\xi))\cup
(P_\ell (\bn^T(p),+\xi)\cap P_\ell(\bn^T(p),-\xi)\right)
,
\]
which is called the {\it spacelike region} with respect to $(\bn^T(p),\xi).$
Here, ${\rm Int} X$ is the interior of $X.$
We also consider the following subsets of $S(\bn^T(p),\xi)$:
\begin{eqnarray*}
S^+(\bn^T(p),\xi)&=&\{ \bx| \langle \bx-f(p),\bv^+_p\rangle\geq 0, \langle \bx-f(p),\bv^-_p\rangle\leq 0\ \mbox{and}\ \langle \bx-f(p),\xi\rangle\geq 0\}.
\end{eqnarray*}
We remark that $\xi \in S^+(\bn^T(p),\xi).$
Then we have the following lemma.
\begin{Lem}
Let $f:M\lon \R^{n+1}_1$ be a spacelike embedding of a closed orientable manifold with $\dim M<n-1.$
If $f$ is L-tight, then
there exists a spacelike affine subspace $V\subset \R^{n+1}_1$ with
$\dim V=n-1$ such that $f(M)\subset V.$
\end{Lem}
\demo
Since $f$ is L-tight, the lightlike tangent hyperplanes  at any point $p\in M$ are
the support plane of $f(M).$
\par
Suppose that there exists $p\in M$ such that
\[
f(M)\subset F_\ell (\bn^T(p),+\xi)\cap F_\ell (\bn^T(p),-\xi),
\]
for any $\xi\in N_1(M)_p[\bn^T].$ 
We arbitrary choose $\xi\in N_1(M)_p[\bn^T].$
Since $HP(\bv^\pm_p,c^\pm)$ are the tangent lightlike hyperplanes at $p\in M,$
we have
$T_{f(p)}f(M)\subset T_SM[\bn^T,\xi]_p.$
By the fact $\dim T_SM[\bn^T,\xi]_p=n-1$ and the assumption $\dim M<n-1,$
there exists $\xi'\in N(M)_p[\bn^T]$ such that $\xi\not= \xi'$ and 
\[
f(M)\subset (F_\ell (\bn^T(p),+\xi')\cap F_\ell (\bn^T(p),-\xi'))\bigcap  (F_\ell (\bn^T(p),+\xi)\cap F_\ell (\bn^T(p),-\xi)).
\]
Therefore, we have $T_{f(p)}f(M)\subset T_SM[\bn^T,\xi]_p\cap T_SM[\bn^T,\xi']_p.$
We can inductively proceed this process, so that we have
\[
 f(p)+T_{f(p)}f(M)\subset \bigcap _{i=1}^\ell T_SM[\bn^T,\xi]_p.
 \]
 However, there exists $\ell $ such that $\dim \bigcap _{i=1}^\ell T_SM[\bn^T,\xi]_p<\dim M.$
 This is a contradiction.
 Therefore, $f(M)\subset S(\bn^T(p),\xi)$ at any point $p\in M.$
 \par
 Suppose that $f(M)\subset S^+(\bn^T(p),\xi)$ at a point $p\in M.$
 Since $\dim M < n-1,$  there exists a closed loop $\bgamma :[0,1]\lon N_1(M)[\bn^T]_p$
 such that $\bgamma(0)=\bgamma (1)=\xi$ and $\bgamma (1/2)=-\xi.$
 By the assumption that $f$ is L-tight, there exists $\bar{\xi}\in N_1(M)[\bn^T]_p$ 
 such that
 \[
 f(M)\subset S^+(\bn^T(p),\bar{\xi})\cap S^+(\bn^T(p),-\bar{\xi})=T_SM[\bn^T,\xi]_p.
 \]
 Here $T_SM[\bn^T,\xi]_p$ is a spacelike affine subspace in $\R^{n+1}_1.$
\enD
Then we have the following theorem.
\begin{Th}Let $f:\mathbb{S}^{s}\lon \R^{n+1}_1$ be a spacelike embedding with $n-1>s$.
Then the following conditions are equivalent:
\par
{\rm (1)} $\tau _\ell (\mathbb{S}^{s},f)=2$,
\par
{\rm (2)} $f$ is L-tight,
\par
{\rm (3)} There exists a spacelike affine subspace $V\subset \R^{n+1}_1$ with $\dim V=s+1$
such that $f(\mathbb{S}^s)$ is a convex hypersurface in $V.$
\end{Th}
\demo
By  Theorem 9.2, the conditions (1) and (2) are equivalent.
It is trivial that the condition (3) implies the condition (2).
We now assume that $f$ is L-tight.
By Lemma 9.3, there exists a spacelike affine subspace $V$ in $\R^{n+1}_1$ with $\dim V=n-1$ such that 
$f(\mathbb{S}^s)\subset V.$
For any $p\in \mathbb{S}^s$ and $\xi\in N_1(\mathbb{S}^s)[\bn^T(p)],$ 
$HP(\bv^\pm,c^\pm)\cap V=V$ or $HP(\bv^\pm,c^\pm)\cap V$ is a hyperplane in $V.$
Since $f$ is L-tight, every tangent hyperplane in $V$ is a support plane of $f(\mathbb{S}^s)$ in $V.$
Therefore, $f(\mathbb{S}^s)$ is tight in $V$ in the Euclidean sense.
Then we can apply the result of submanifolds in the Euclidean space \cite{CL57}, so that there exists a a spacelike affine subspace $V\subset \R^{n+1}_1$ with $\dim V=s+1$
such that $f(\mathbb{S}^s)$ is a convex hypersurface in $V$.
This completes the proof.
\enD
\section{Special cases}
In this section we consider submanifolds in Euclidean space and Hyperbolic space as
special cases as the previous results.
\subsection{Submanifolds in Euclidean space}
Let $\R^n_0$ be the Euclidean space which is given by the equation $x_0=0$ for $\bx =(x_0,x_1,\dots ,x_n)\in\R^{n+1}_1$.
Consider an immersion $f:M\lon \R^n_0$, where $M$ is a closed orientable manifold.
We remember that the total absolute lightcone curvature $K^*_{\ell}(p)$ is independent of the choice of $\bn^T$ (cf., \S 6).
Therefore, we can adopt $\bn^T=\be _0=(1,0,\dots ,0)$ as a future directed timelike unit normal
vector field along $f(M)$ in $\R^{n+1}_1.$ 
In this case $N_1(M)[\bn^T]=N_1(M)[\be_0]$ is the unit normal bundle $N^e_1(M)$ of $f(M)$ in $\R^n_0$
in the Euclidean sense.
Therefore, the lightcone Gauss map $\widetilde{\mathbb{LG}}(\bn^T)$ is given by
$\widetilde{\mathbb{LG}}(\bn^T)(p,\bxi)=\be_0+\bxi=\be_0+\mathbb{G}(p,\bxi),$
where $\mathbb{G}:N^e_1(M)\lon \mathbb{S}^{n-1}$ is the Gauss map of the unit normal bundle $N^e_1(M)$\cite{CL57}.
Since $\be_0$ is a constant vector, we have
\[
K^*_{\ell}(p)=K^*(p), 
\]
where $K^*(p)$ is the total absolute curvature of $M$ at $p$ (cf., \cite{CL57}) in the Euclidean sense.
Therefore, Theorem 7.3 is the original Chern-Lashof theorem in \cite{CL57}.
If $\dim M=n-1$, then the $\widetilde{K}_\ell^\pm(p)=(\pm 1)^{n-1} K(p)$ where $K$ is the Gauss-Kronecker curvature of $M.$ Thus, if $n$ is odd, then $\widetilde{K}_\ell^\pm(p)=K(p)$.
Moreover, $|\widetilde{K}_\ell^\pm|(p)=|K|(p)$ for any $n.$
Therefore, Theorems 8.1, 8.3 and 8.4 are the original integral formulae in the Euclidean sense\cite{CL57}.
Furthermore, if $n=3$, then the Proposition 8.6 is the original Willmore inequality in Euclidean space 
\cite[Theorem 7.2.2]{Will}.
\par
On the other hand, the intersection of a lightlike hyperplane with $\R^n_0$ is a hyperplane in
$\R^n_0$, so that the notion of lightlike tightness is equivalent to the original notion of 
the tightness\cite{CR85}.
\par
We remark that if $\bn^T=\bv$ is a constant timelike unit vector, the spacelike submanifold $f(M)$ is a
submanifold in the spacelike hyperplane $HP(\bv,c)$. Since $HP(\bv,c)$ is isometric to the Euclidean space
$\R^n_0,$ all results for the case $\bn=\be _0$ hold in this case.

\subsection{Submanifolds in Hyperbolic space}
Let $f:M\lon \mathbb{H}^n(-1)$ be an immersion into the hyperbolic space. Then we adopt $\bn^T(p)=f(p)$.
In this case $N_1(M)[\bn^T]=N_1(M)[f]$ is the unit normal bundle $N^h_1(M)$ of $f(M)$ in $\mathbb{H}^n(-1).$
Therefore, the lightcone Gauss map $\widetilde{\mathbb{LG}}(\bn^T)$ is given by
$\widetilde{\mathbb{LG}}(\bn^T)(p,\bxi)=\widetilde{f(p)+\bxi}=\widetilde{\mathbb{L}}(p,\bxi),$
where $\widetilde{\mathbb{L}}:N^h_1(M)\lon \mathbb{S}^{n-1}_+$ is the horospherical Gauss map of the unit normal bundle $N^h_1(M)$(cf., \cite{BIR10}).
Thus, we have
\[
K^*_{\ell}(p)=\int _{N^h_1(M)_p} |\widetilde{K}_h(p,\bxi)|d\sigma _{k-2},
\]
which is the total absolute horospherical curvature of $M$ at $p$ (cf., \cite{BIR10}) in $\mathbb{H}^n(-1).$
Therefore, $\tau _\ell (M,f)=\tau _h(f)$.
\par
On the other hand,
let $f:M\lon \mathbb{H}^n(-1)$ be an embedding such that $M$ is a closed orientable manifold with ${\rm dim}\, M=n-1.$
In this case, $f(M)$ is a spacelike submanifold of codimension two in $\R^{n+1}_1,$ then
we have
$\tau _\ell ^\pm (M,f)=\tau _h^\pm (f;M)$ (cf., \cite{BIR11}).
In \cite{BIR11} we gave an example of a curve in $\mathbb{H}^2(-1)$ such that
$\tau _h^+(f;M)\not= \tau _h^-(f;M).$
This example can be easily generalized into any higher dimensional case.
\par
On the other hand, the notion of the lightlike tightness is equivalent to the notion of the 
horo-tightness in $\mathbb{H}^n(-1)$\cite{BIR11,CR79,CR85,SolTeu}.
Since the intersection of $\mathbb{H}^n(-1)$ with a spacelike affine subspace $V$ is a round hypersphere in $V,$
the condition (3) in Theorem 9.4 can be changed into the following condition:
\par
(3$'$) $f(\mathbb{S}^s)$ is a metric (round) sphere in $\mathbb{H}^n(-1).$
\par\noindent
Therefore, Theorems 9.2 and  9.4 are characterizations of the horo-tight spheres in $\mathbb{H}^n(-1)$ \cite{BIR11}.
Further results on horo-tight immersions into $\mathbb{H}^n(-1)$ are presented in \cite{SolTeu}.

\begin{flushleft}
\textsc{
Shyuichi Izumiya
\\ Department of Mathematics
\\ Hokkaido University
\\ Sapporo 060-0810
\\ Japan} 
\\
{E-mail address}: izumiya@math.sci.hokudai.ac.jp

\end{flushleft}

\end{document}